\documentclass{amsart}
\usepackage{amsmath, rlepsf,amssymb, amsfonts, mathrsfs}

\theoremstyle{plain}
\newtheorem{thm}{Theorem}[section]
\newtheorem{lem}[thm]{Lemma}
\newtheorem{cor}[thm]{Corollary}

\newtheorem{rem}[thm]{Remark}

\newtheorem{fct}[thm]{Fact}

\newtheoremstyle{definition}{7pt plus6.3pt minus6.3pt}{7pt plus3pt minus3pt}%
{\rm}{}{\bf}{}{0.75em}{\thmname{#1}\thmnumber{ #2}\thmnote{\sl\stdspace#3}}
\theoremstyle{definition}
\newtheorem{example}[thm]{Example}
\newtheorem{exercise}[thm]{Exercise}

\newcommand{\bbr}{\begin{rem}\em} 
\newcommand{\eer}{\end{rem}}

\newcommand{\bex}{\begin{example}} 
\newcommand{\eex}{\end{example}}

\newcommand{\bhw}{\begin{exercise}} 
\newcommand{\ehw}{\end{exercise}}

\newcommand{\be}{\begin{enumerate}}
\newcommand{\ee}{\end{enumerate}}

\def\C{\hbox{$\mathbb C$} }
\def\Z{\hbox{$\mathbb Z$} }

\def\R{\hbox{$\mathbb R$} }

\def\co{\colon\thinspace}

\def\dfn#1{{\em #1}}

\begin{document}


%

\title{Lectures on Contact Geometry in Low-Dimensional Topology}
\author{John B. Etnyre}

\address{Georgia Institute of Technology, Atlanta, GA 30332-0160}
\email{etnyre@math.gatech.edu}
\urladdr{http://www.math.gatech.edu/\char126 etnyre}
\thanks{The author thanks the PCMI organizers for the opportunity to give these lectures and David Shea Vick and Elena Bogdan for providing valuable comments on the first draft of the notes. This work was supported in part by NSF CAREER Grant (DMS--0239600) and FRG-0244663. }



\begin{abstract}
This article sketches various ideas in contact geometry that have become useful in low-dimensional topology. Specifically we (1) 
outline the proof of Eliashberg and Thurston's results concerning perturbations
of foliatoins into contact structures, (2) discuss Eliashberg and Weinstein's symplectic handle attachments, and (3) briefly discuss Giroux's insights into 
open book decompositions and contact geometry. Bringing these pieces together we discuss the
construction of  ``symplectic caps'' which are a key tool in the application of contact/symplectic
geometry to low-dimensional topology. 
\end{abstract}

\maketitle


\section{Introduction}
Contact geometry has been a key tool in many recent advances in low-dimensional topology. For example contact geometry was an integral part
in the following results:
\begin{enumerate}
\item Kronheimer and Mrowka's proof that all non-trivial knots satisfy property $P$ \cite{KronheimerMrowka04}.\\
(Recall, a knot satisfies property $P$ if non-trivial surgery on it yields a manifold with non-trivial fundamental group.)
\item Ozsv\'ath and Szab\'o's proof that the unknot, trefoil and figure eight knots are all determined by surgery \cite{OzsvathSzabo04a, OzsvathSzabo??}. That is, if $r$ Dehn surgery on $K$ is the same
as $r$ surgery on the unknot, trefoil or figure eight knots then the knot was the unknot, trefoil or figure eight knot. This result, for the unknot, was originally proven
by Kronheimer-Mrowka-Ozsv\'ath and Szab\'o \cite{KronheimerMrowkaOzsvathSzabo??}.
\item Ozsv\'ath and Szab\'o's proof that Heegaard-Floer invariants detect the Thurston norm of a manifold and the minimal Seifert genus of a knot \cite{OzsvathSzabo04a}.
\end{enumerate}
There are many other results in which contact geometry has played a key role. For a brief, far from complete list of applications we have:  achiral Lefschetz fibrations \cite{EtnyreFuller06}, Harer's conjecture on fibered knots  \cite{GirouxGoodman06},  characterizing fibered knots in terms of Heegaard-Floer theory \cite{Ghiggini??, Ni??}, new knot invariants \cite{Ng05} and
the existence of Engel structures \cite{Vogel??}. 

We outline how contact geometry shows up in the results highlighted above by sketching the proof that Heegaard-Floer homology detects the Thurston norm
(that is the minimal genus of a surface representing a homology class in a 3-manifold).  Any unfamiliar terminology will be explained in the sections below.

Start with a closed irreducible oriented 3-manfiold $M$ and a surface $\Sigma\subset M$ such that $\Sigma$ is of minimal genus among surfaces homologous to it. (Assume the genus of $\Sigma$ is larger than 0.)
\begin{enumerate}
\item A theorem of Gabai \cite{Gabai83}, Theorem~\ref{GabaiThm} below, gives a taut foliation $\mathcal{F}$ on $M$ that contains $\Sigma$ as a leaf.
\item A theorem of Eliashberg and Thurston \cite{EliashbergThurston98}, Theorem~\ref{thm:f2c} below, gives a positive and negative contact structure $\xi_\pm$ on $M$ that is $C^0$-close
to $\mathcal{F}.$
\item They also show, see the proof of Theorem~\ref{makefill} below,  how to construct a symplectic structure $\omega$ on $X=M\times[-\epsilon, \epsilon]$ that weakly fills 
\[
(M,\xi_+)\cup (-M,\xi_-).
\] 
\item A theorem of Eliashberg \cite{Eliashberg04} and, independently, the author \cite{Etnyre04a}, Theorem~\ref{eet} below, shows how to find a closed symplectic manifold $(X',\omega')$ into
which $(X,\omega)$ embeds.  This closed manifold is constructed by finding ``symplectic caps'' to cap off the boundary of $X.$ See Figure~\ref{Xconstruction}.
\begin{figure}[ht]
  \relabelbox \small {\centerline{\epsfbox{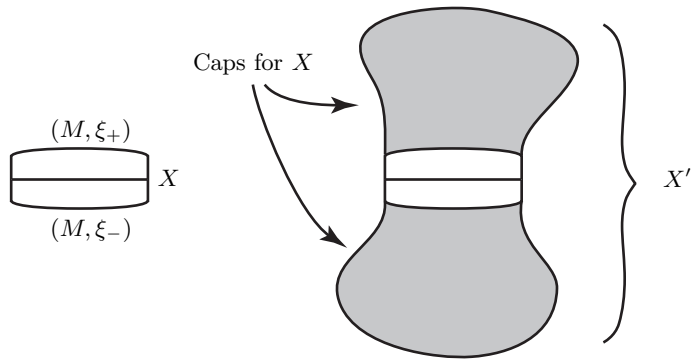}}} 
  \relabel{1}{$X$}
  \relabel{2}{$(M,\xi_+)$}
  \relabel{3}{$(M,\xi_-)$}
  \relabel{4}{Caps for $X$}
  \relabel{5}{$X'$}  
  \endrelabelbox
        \caption{On the left is the symplectic manifold $X=M\times[-\epsilon,\epsilon]$. On the right, is the symplectic manifold $X'$ that is constructed by ``gluing'' caps
        onto $X.$}
        \label{Xconstruction}
\end{figure}
There are various ways to construct these caps. We discuss the constructions in \cite{Etnyre04a} which uses:
\begin{enumerate}
\item Giroux's correspondence between open book decompositions and contact structures \cite{Giroux02}, Theorem~\ref{GirouxCor} below, and
\item Eliashberg \cite{Eliashberg90a} and, independently, Weinstein's \cite{Weinstein91} ideas of contact surgery and symplectic handle attachment, see Theorem~\ref{thm:ew} below.
\end{enumerate}
\item At this point one needs to use something like Seiberg-Witten theory or Heegaard-Floer theory to conclude something about $M$ or $\Sigma$ based on
the existence of $(X',\omega')$ with $(M,\xi_+)$ nicely embedded in it. Recall, there are lots of non-vanishing theorems for symplectic manifolds and with
non-vanising invariants you can frequently conclude things about surfaces in homology classes using an adjunction type inequality. Specifically, in our current
situation we know that the Heegaard-Floer invariant of $X'$ is non-zero in the $\text{spin}^c$ structure associated to the symplectic form and this implies that
the Heegaard-Floer homology $HF^+(M,\frak{s}_\xi)\not=0,$ where $\frak{s}_\xi$ is the $\text{spin}^c$ structure associated to $\xi.$ Recall that for any 
$\text{spin}^c$ structure $\frak{s}$ for which $HF^+(M,\frak{s})\not=0$ we have the adjunction inequality 
\[
|\langle c_1(\frak{s}), [\Sigma]\rangle|\leq 2g-2.
\]
But since $\Sigma$ is a leaf of $\mathcal{F}$ it is easy to see that $\langle c_1(\frak{s}_\xi), [\Sigma]\rangle= 2g-2.$ Thus we see that we can detect the 
minimal genus of a non-trivial homology class in $M$ by seeing how $c_1$ of all the $\text{spin}^c$ structures on $M$ with non-zero Heegaard-Floer homologies 
evaluate on the homology class. In other words, Heegaard-Floer homology detects the Thurston norm of homology classes in $M.$
\end{enumerate}
It is not too hard to adapt this line of argument to see that Heegaard-Floer homology also detects the Seifert genus of a knot in $S^3$. The main idea is to
perform 0 surgery on the knot and then apply the outline above to the surgered manifold. The other results mentioned above require somewhat different 
arguments, but the difference is in step (5) and hence the contact geometric input is largely the same. 

These lectures are devoted to understanding the contact geometric part of the above outline, that is steps (2) through (4). We begin in Section~\ref{terms}
by defining all our basic concepts like contact structure, foliation, tight, symplectically fillable and so on. In Section~\ref{sec:f2c} we give a detailed sketch of the
proof of Eliashberg and Thurston's theorem (step (2) in the above outline). In Section~\ref{taut2fill} we discuss the relation between various notions associated to 
foliations and similar notions associated to contact structures. This discussion ends by constructing the symplectic manifold $(X,\omega)$ from step (3) of the
above outline. The next two sections are devoted to the construction of symplectic caps in step (4). Specifically, we discuss Legendrian surgery and symplectic
handle attachment in Section~\ref{shandles}. The following section states Giroux's correspondence between contact structures and open book decompositions.
We then use this correspondence to outline the construction of symplectic caps. The material in this last section is a bit more sketchy than that in the previous sections.
This is because there are already lecture notes devoted to this topic. See \cite{EtnyreOBN} for a more complete discussion of Giroux's correspondence and the construction
of symplectic caps. Thus you can regard these lectures as primarily concerning the relation between contact structures an foliations as exemplified in steps (2) and (3) above. Step (1) is purely foliation theoretic and would take us far afield of contact geometry. Fore more details on step (1) see Gabai's original work
\cite{Gabai83} or the text book \cite{CandelConlon03}. Step (5) is also beyond the scope of these lectures, in that it does not explicitly use contact geometry. For more on this step we 
refer the reader to the original research articles or the more expository articles \cite{Geiges??, OzsvathSzabo????}. 
We would like to emphasize that while we are using the applications to low-dimensional topology to motivate and guide our discussion of
various topics in contact geometry, each of these topics is quite beautiful in its own right and there are many exciting directions for future research in all of these
topics. We hope this beauty and potential is evident in what follows!

\section{Contact structures and foliations}\label{terms}
\begin{figure}[ht]
  \relabelbox \small {\centerline{\epsfbox{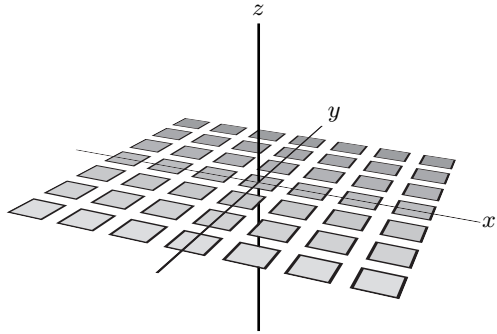}}} 
  \relabel{1}{$x$}
  \relabel{2}{$y$}
  \relabel{3}{$z$}
  \endrelabelbox
        \caption{The plane field $\xi_1=\ker dz.$}
        \label{fig:basicexfol}
\end{figure}

Throughout this section $M$ will be a closed oriented 3--manifold. A \dfn{plane field} on $M$ is simply a 2--dimensional
sub-bundle of the tangent bundle of $M.$ That is, at each $x\in M$ we have a plane $\xi_x$ in the tangent bundle $T_xM.$ 
Locally, one can always find a 1--form $\alpha$ such that 
\[
\xi_x=\ker \alpha_x,
\]
for each $x\in M.$
\bhw
Show that $\alpha$ may be chosen to be a global 2--form if and only if $\xi$ is orientable. \\
Hint: It might be helpful to use a Riemannian metric on $M.$
\ehw
\bex
Consider $\R^3$ with coordinates $(x,y,z).$ Let $\alpha_1=dz$ and set 
\[\xi_1=\ker \alpha_1=\text{span}\{\frac{\partial}{\partial x}, \frac{\partial}{\partial y}\}.\] 
See Figure~\ref{fig:basicexfol}.

\begin{figure}[ht]
  \relabelbox \small {\centerline{\epsfbox{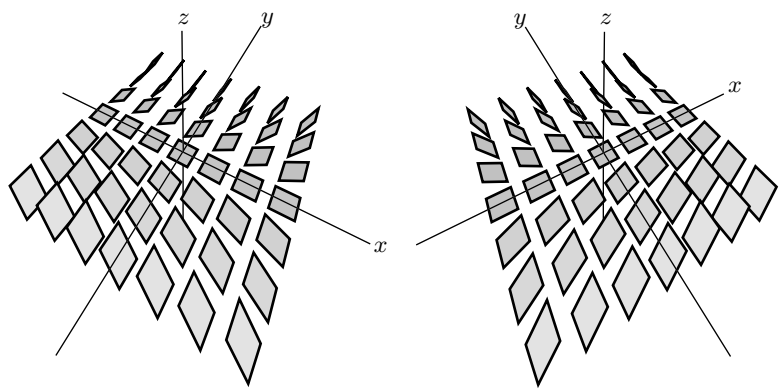}}} 
  \relabel{1}{$x$}
  \relabel{2}{$y$}
  \relabel{3}{$z$}
  \relabel{4}{$x$}
  \relabel{5}{$y$} 
  \relabel{6}{$z$} 
  \endrelabelbox
        \caption{On the left is $\xi_2=\ker (dz -y\, dx)$ and on the right is $\xi_3=\ker(dz+y\, dx).$}
        \label{fig:basicex}
\end{figure}

This is a fairly simple plane field. It is constant and does not change from point to point.
\eex
\bex\label{mainex}
Again consider $\R^3$ with coordinates $(x,y,z).$ Let $\alpha_2=dz-y\, dx$ and $\alpha_3=dz+y\, dx.$ Set 
\[\xi_2=\ker \alpha_2=\text{span}\{ \frac{\partial}{\partial x}+y\frac{\partial}{\partial z}, \frac{\partial}{\partial y}\} \]
and 
\[\xi_3=\ker \alpha_3=\text{span}\{ \frac{\partial}{\partial x}-y\frac{\partial}{\partial z}, \frac{\partial}{\partial y}\}.\] 
See Figure~\ref{fig:basicex}.
So each $\xi_2$ and $\xi_3$ consist of horizontal planes (that is, parallel to the $xy$-plane) 
at any point in the $xz$-plane and as you leave the $xz$-plane along a
ray perpendicular to the $xz$-plane the planes $\xi_2$ and $\xi_3$ are always tangent to this ray and twisting a total of 
$90^\circ$ in a counterclockwise, respectively clockwise, manner.
\eex

We call a plane field $\xi$ a \dfn{foliation} if there is a 1--form $\alpha$ with $\xi=\ker \alpha$ for which
\[
\alpha\wedge d\alpha=0.
\]
We call $\xi$ a \dfn{positive contact structure}, respectively \dfn{negative contact structure}, on $M$ if there is a  1--form
$\alpha$ with $\xi=\ker \alpha$ for which
\[
\alpha\wedge d\alpha \text{ is a positive multiple of the volume form on $M$},
\]
respectively
\[
\alpha\wedge d\alpha \text{ is a negative multiple of the volume form on $M$}.
\]
These conditions are usually stated
\[
\alpha\wedge d\alpha>0
\]
and
\[
\alpha\wedge d\alpha<0.
\]
A plane field $\xi$ is called a \dfn{positive confoliation}, respectively \dfn{negative confoliation}, if there is a 1--form $\alpha$ with
$\xi=\ker \alpha$  for which
\[
\alpha\wedge d\alpha\geq 0,
\]
respectively 
\[
\alpha\wedge d\alpha\leq 0.
\]
\bhw
Show that these definitions do not depend on the 1--form $\alpha$ chosen to define $\xi.$
\ehw
\bhw
Show that if a 3--manifold supports a contact structure then it is orientable.
\ehw
\begin{thm}[Frobenius, see \cite{Conlon01}]
If  a plane field $\xi$ is closed under Lie brackets (that is, if $v, w$ are sections of $\xi$ then their Lie bracket $[v,w]$ is
also a section of $\xi$), then $M$ is foliated by surfaces tangent to $\xi.$ 
\end{thm}
To make the theorem more precise we say a 3-manifold $M$ is \dfn{foliated by surfaces} if $M$ can be written as the disjoint union of surfaces such
that each point on $M$ is contained in a coordinate chart that maps the intersection of the surfaces with the chart  
to the constant $z$-hyperplanes in $\R^3$ (with
coordinates $(x,y,z)$). The degree of smoothness of the foliation is the degree of smoothness with which 
the transition functions between charts can be chosen. 
The surfaces that show up in this theorem are called \dfn{leaves} of the foliation. The Frobenius theorem is an ``integrability'' theorem.
That is you can ``integrate'' information in the tangent space into the manifold: given $\xi$ in the tangent space you find things, surfaces,
that are actually in the manifold that induce the plane field. 
\bhw
Show that the condition $\alpha\wedge d\alpha=0$ is equivalent to $\xi=\ker \alpha$ being closed under Lie bracket.\\
Hint: Recall the formula for $d\alpha$ in terms of $\alpha$ and Lie brackets.
\ehw
From this exercise we see that the 1--form conditions in the definitions of foliations and contact structures have to do with the plane field being tangent
to surfaces. This integrability condition implies that a contact structure cannot be tangent to any surface along an open set of the surface.
(Of course, a surface can be tangent to a contact structure at isolated points and even along lines in the surface, but not along an open
set in the surface.)
\bex
Back to the example $\R^3$ with $\xi_1=\ker \alpha_1,$ where $\alpha_1=dz.$ Clearly $d\alpha_1=d(dz)=0$ so $\alpha_1\wedge d\alpha_1
=0$ and $\xi_1$ is a foliation. If we define the surfaces $S_{z_0}=\{(x,y,z)| z=z_0\}$ in $\R^3,$ then clearly
\[
\R^3=\coprod_{z\in \R} S_{z},
\]
\begin{figure}[ht]
  \relabelbox \small {\centerline{\epsfbox{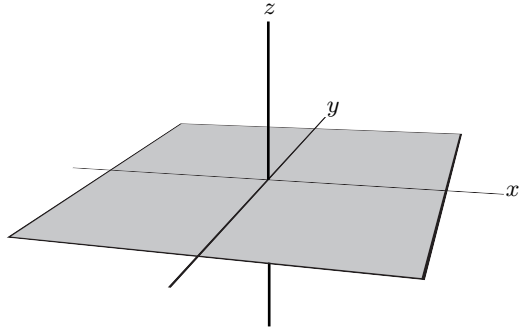}}} 
  \relabel{1}{$x$}
  \relabel{2}{$y$}
  \relabel{3}{$z$}
  \endrelabelbox
        \caption{Here is the surface to which the planes in Figure~\ref{fig:basicexfol}  are tangent.}
        \label{fig:basicfol}
\end{figure}
and $T_{(x,y,z)}S_z=(\xi_1)_{(x,y,z)}.$ See Figure~\ref{fig:basicfol}. This verifies the Frobenius theorem in this example.
\eex
\bex
Let $M=S^1\times \Sigma$ for some surface $\Sigma.$ Let $\theta$ be the coordinate on the $S^1$ factor and $\alpha=d\theta.$ Clearly the
plane field $\xi=\ker\alpha$ is always tangent to $\{pt\}\times \Sigma$ so $\xi$ is a foliation with leaves the surface fibers. More generally, we could let
$M$ be a $\Sigma$ bundle over $S^1$ and again let $\xi$ be the tangents to the fibers. Note $d\alpha=d(d\theta)=0$ so $\alpha\wedge d\alpha=0$ as we
expect for a foliation.
\eex
\bex
The plane fields $\xi_2$ and $\xi_3$ from Example~\ref{mainex} are positive and, respectively, negative contact structures since one can easily check that
$\alpha_2\wedge d\alpha_2>0$ and $\alpha_3\wedge d\alpha_3<0.$
\eex
\bex
We now give a contact structure on a closed 3-manifold. Let $S^3$ be the unit 3-sphere in $\R^4=\C^2.$ Set $\alpha=r_1^2d\theta_2+ r_2^2d\theta_2$ where
the coordinates on $\C^2$ are $(z_1,z_2)$ and $z_j=r_je^{i\theta_j}.$ 
\bhw
Check that when $\alpha$ is restricted to $S^3$ then $\alpha\wedge d\alpha>0.$ Also show that $\ker\alpha$ is the set of complex tangencies to $S^3,$ that is 
$\xi$ is the set of vectors $v$ tangent to $S^3$ for which $iv$ is also tangent to $S^3,$ where $i=\sqrt{-1}.$
\ehw
\eex

\begin{fct}
All oriented 3--manifolds have foliations and positive (negative) contact structures.
\end{fct}
There are many proofs of this fact. We will indicate one proof in Section~\ref{sec:ob}. For another approach see \cite{Etnyre03}.

We would now like to understand plane fields locally. For this we start with the following lemma.
\begin{lem}\label{lem:locform}
Given a plane field $\xi$ there are local coordinates so that a 1--form $\alpha$ with $\xi=\ker \alpha$ has the
form
\[
\alpha= dz- a(x,y,z)\, dx
\]
\end{lem}

\bhw
Prove this lemma.\\
Hint: Let $\phi\co [-1,1]\times [-1,1] \to M$ be an embedding of a disk such that $\phi_*\frac{\partial}{\partial z}$ is transverse to
$\xi,$ where the coordinates on the disk are $(x,z).$ Now find a vector field $v$ in $\xi$ that is transverse to the image of this
embedding. Use the vector field to extend $\phi$ to an embedding of a 3--ball.
\ehw
\bbr
Note that your proof (if you followed the hint)
of Lemma~\ref{lem:locform} can be extended to the following setting: if $\gamma\co N\to M$
is any embedding, where $N=[0,1]$ or $S^1,$ 
for which $\gamma'(x)\in \xi_{\gamma(x)}$ for all $x\in N$ then we can extend $\gamma$ to $\Gamma\co
N\times[-\epsilon,\epsilon]\times[-\epsilon,\epsilon]\to M$ so that $\alpha$ has the desired form in these coordinates. 
Also note that if $N=[0,1]$ then you can also assume the coordinate on $\gamma$ is $y$ and extend $\gamma$ to $\Gamma\co 
[-\epsilon,\epsilon]\times N\times[-\epsilon,\epsilon]\to M$ so that $\alpha$ has the desired form in these coordinates. 
\eer

To see how this local form for $\alpha$ relates to $\xi$ being a foliation or contact structure, we have the following lemma.
\begin{lem}\label{lem:type}
Given a plane field $\xi$ and a 1--form $\alpha$ as in Lemma~\ref{lem:locform}, then 
\begin{enumerate}
\item $\xi$ is a positive (negative) contact structure if and only if
\[
\frac{\partial a}{\partial y}>0 \quad\left(\frac{\partial a}{\partial y}<0\right).
\]
\item $\xi$ is a foliation if and only if 
\[
\frac{\partial a}{\partial y}=0.
\]
\end{enumerate}
\end{lem}
\bhw
Verify this lemma.
\ehw
One can improve the above local results for the 1--form $\alpha$ when $\ker \alpha$ is a contact structure or foliation. 
\begin{thm}[Darboux/Pfaff Theorem, see \cite{Conlon01, McDuffSalamon98}]
If $\xi$ is a foliation then there are local coordinates $(x,y,z)$ such that
\[
\xi=\ker dz.
\]
If $\xi$ is a positive (negative) contact structure then there are local coordinates $(x,y,z)$ such that
\[
\xi=\ker (dz - y\, dx),\quad \left(\xi=\ker(dz + y\, dx)\right).
\]
\end{thm}
This theorem says that ``locally all foliations (positive/negative contact structures) look the same''. Note this is very different form other 
``geometries'', like Riemannian metrics, which can look very different locally. Thus foliations and contact structures are in some sense  
insensitive to local things and thus if they tell us anything about the manifold, it will have to be something global.

This theorem also indicates similarities between foliations and contact structures ({\em i.e.} they both have local normal forms). We will 
see many more similarities below, but there are some differences. One major difference is that there are no ``non-tivial deformations'' of a
contact structure where foliations have ``non-tirival deformations''. To better understand this last sentence consider the following theorem.
\begin{thm}[Gray's Theorem, \cite{McDuffSalamon98}]
If $\xi_t, t\in [0,1],$ is a 1-parameter family of contact structures on $M$ that agree off of a compact subset of $M$  then
there is a 1-parameter family of diffeomorphisms $\psi_t\co  M \to M$ such that $(\psi_t)_*\xi_0=\xi_t.$
\end{thm}
This theorem says that isotopies of contact structures (as plane fields) are equivalent to isotopies of the manifold. 


If $\xi_1$ and $\xi_2$ are two contact structures on $M$ then a diffeomorphism $f\co M\to M$ is called a \dfn{contactomorphism}
from $\xi_1$ to $\xi_2$ if $f_*(\xi_1)=\xi_2.$ So Grey's Theorem says that any family of contact structures are related by a 
family of contactomorphisms. The situation for foliations is quite different.

\bex
Let $\mathcal{F}_s$ be the foliation of $T^2$ by lines of slope $s.$ Let $\xi_s=\mathcal{F}_s\times S^1$ be the product foliation
on $T^3=T^2\times S^1.$ 
\eex
\bhw
Show there is no family of diffeomorphisms $\phi_s\co T^3\to T^3.$ such that $(\phi_s)_*\xi_0=\xi_s.$
\ehw
This last example (and exercise) show that there are deformations of foliations that do not come from diffeomorphisms of the underlying
3--manifold. 

To see further similarities between foliations and contact structures we consider special foliations and contact structures 
on $D^2\times S^1.$
\bex
In this example we construct a Reeb foliation on the solid torus. We will construct a foliation on $\R^3$ that is invariant under translations in the
$z$-direction. Thus we can look at $\R^2\times S^1$ thought of as $\R^3$ modulo $z\mapsto z+1.$ Then the unit disk in the $xy$-plane times $S^1$ will
be a solid torus. Consider the function
\[
f(x,y,z)= g(x^2+y^2)e^z
\]
where $g\colon [0,\infty) \to \R$ is a strictly decreasing function equal to 1 at 0 and 0 at 1 (moreover assume all the higher derivatives of $g$ are 0 at 0). One may 
check that $f$ has no critical points so the level sets of $f$ give a foliation of  $\R^3.$ Moreover it is easy to check that the foliation is invariant under translation
in the $z$-direction. Thus we get a foliation induced on the unit disk in the $xy$-plane cross $S^1.$
The unit circle in the $xy$-plane times $S^1$ will be a closed torus leaf in the foliation. The interior 
of the solid torus is foliated by leaves diffeomorphic to $\R^2.$
\bhw
Picture this foliation.
\ehw

\begin{figure}[ht]
  \relabelbox \small {\centerline{\epsfbox{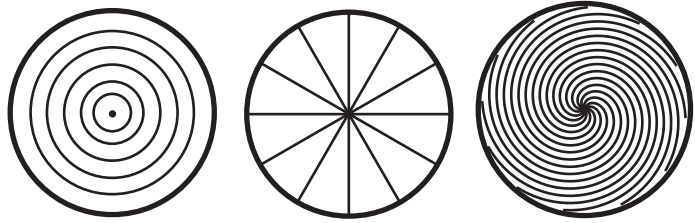}}} 
    \endrelabelbox
        \caption{The meridional disk in the Reeb torus is shown on the left. The meridional disk in the Lutz tube is shown in the middle and on the right is the
        meridional disk in the Lutz tube with its interior pushed up slightly.}
        \label{disks}
\end{figure}

If we let $D=D^2\times\{pt\}$ be a meridional disk in the solid torus then intersecting $D$ with $\xi$ will induce a
singular foliation on $D.$ This foliation is shown on the left of Figure~\ref{disks}.
\eex
A foliation on a manifold $M$ is said to have \dfn{Reeb components} if there is a solid torus in $M$ such that the
foliation is homeomorphic with the one in the pervious example. (Note we are thinking of a foliation here in terms of its leaves and not
as a plane field. In this way it makes sense to talk about homeomorphic foliations.)
A foliation on $M$ that does not have any Reeb components is said to be \dfn{Reebless}.

\begin{figure}[ht]
  \relabelbox \small {\centerline{\epsfbox{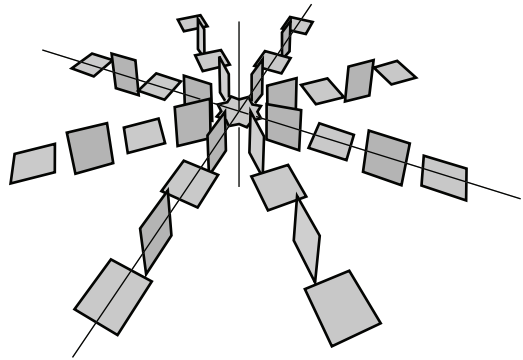}}} 
    \endrelabelbox
        \caption{The contact structure $\xi_{ot}$ on $\R^3.$}
        \label{fig:std_ot}
\end{figure}

Before we move to the next example let's formalize this idea of an ``induced foliation''. Let $\xi$ be any plane field on
a 3--manifold $M.$ If $\Sigma$ is a surface embedded in $M$ then at each point $x\in \Sigma,$ $\xi_x\cap T_x\Sigma$
is either a line in $T_x\Sigma$ or all of $T_x\Sigma.$ The points where $\xi_x=T_x\Sigma$ are called \dfn{singular
points}. Away form the singular points we have a line field on $\Sigma.$ It is easy to use the Forbenious Theorem (or 
in this case just the existence of solutions to ordinary differential equations) to see that we can foliate $\Sigma$ with
1-manifolds away from the singular points. This is called the \dfn{induced singular foliation}, or when $\xi$ is a
contact structure it is sometimes called the \dfn{characteristic foliation}.
\bex
In this example we construct a Lutz tube, this is a special contact structure on a solid torus. Again we will do this by constructing a contact structure on $\R^3$ that is
invariant under translations in the $z$-direction. Then if we consider $\R^3/\sim,$ where $(x,y,z)\sim (x,y,z+1),$ we get a contact structure on $\R^2\times S^1.$
Choosing a disk in $\R^2$ will give our desired contact structure on a solid torus. Let $\alpha=\cos r\, dz + r\sin r\, d\theta.$ This is a 1-form in cylindrical coordinates.
Set
\[
\xi_{ot}=\ker \alpha.
\] 
This contact structure is shown in Figure~\ref{fig:std_ot}. 
Note $\xi_{ot}$ is radially symmetric and invariant under translation in the $z$-direction. Also note that along a ray perpendicular to the $z$-axis the contact planes twist in a left handed fashion, and they twist infinitely often. 
The torus $T=\{(r,\theta,z)| r\leq \pi\}$ in $\R^3/\sim$ with the contact structure induced from $\xi_{ot}$ is called a \dfn{Lutz tube}.

If $D=D^2\times \{pt\}$ is a meridional disk in a Lutz tube then the characteristic foliation induced on $D$ by 
the contact structure $\xi$ is shown in the middle of Figure~\ref{disks}. This characteristic foliation is a bit strange in that the boundary of the disk
consists entirely of singularities. This is a very ``non-generic'' phenomenon, that is, usually you don't see lines of 
singularities. In particular, try the following exercise. 
\bhw
If the interior of $D$ is pushed up a little then show the characteristic foliation is as shown on the right hand side of 
Figure~\ref{disks}. By ``pushed up a little'' we mean let $D'=\{(r,\theta,z)| r\leq \pi, z= \epsilon(-r^2+\pi^2)\},$ where
$\epsilon$ is some positive number very very close to zero.
Note there is only one singularity in the new characteristic foliation.
\ehw
\eex
The disk shown on the right hand side of Figure~\ref{disks} is called an \dfn{overtwisted disk}. The key feature of this disk is that its boundary is 
tangent to $\xi$ but the disk itself is transverse to $\xi$ along the boundary. Any contact
structure that contains an overtwisted disk is called \dfn{overtwisted}, otherwise it is called \dfn{tight}.
\bbr
Any overtwisted contact structure contains a Lutz tube (this is not obvious!). So we could have defined
a tight contact structure to be a ``Lutzless'' contact structure in analogy with a Reebless foliation. 
\eer
It is very easy to construct contact structures with Lutz tubes and foliations with Reeb components. In particular, 
for contact structures we have the following result. 
\begin{thm}[Eliashberg 1989, \cite{Eliashberg89}]
On a closed oriented 3--manifold there is a one-to-one correspondence between homotopy classes of plane 
fields and isotopy classes of overtwisted contact structures. 
\end{thm}
\bhw
Using Eliashberg's theorem show that any closed oriented 3--manifold has infinitely many distinct overtwisted 
contact structures.
\ehw
\bhw
Try to show that any overtwisted contact structure contains a Lutz tube.\\
Hint: Try to fine a way of modifying a contact structure to introduce a Lutz tube and then try to see if you
can do this so as not to change the homotopy class of plain field. 
\ehw

So overtwisted contact structures are fairly well understood and seem to be very flexible. It turns out tight contact
structures are much more interesting, but do they exist? The answer is ``yes'' but not always. 
\begin{thm}[Etnyre-Honda, \cite{EtnyreHonda01a}]
Let $M$ be the Poincar\'e homology sphere with its non-standard orientation. (In other words, $M$ can be described as 1 surgery on 
the right handed trefoil knot.) Then $M$ does not admit a 
positive tight contact structure. Moreover, $M\# (-M)$ does not admit any tight contact structures (positive or
negative). Here $-M$ means $M$ with the opposite orientation. 
\end{thm}
This is a little worrisome, but we will see shortly that there are lots of tight contact structures; however, let's first observe an
important property of them that indicates they see subtle properties of topology. For this recall that an oriented 2-dimensional  bundle, like $\xi,$ has an Euler class $e(\xi)\in H^2(M;\Z).$ Moreover, if $\Sigma$ is not closed but transverse to $\xi$ then
choose a vector field $v_x\in \xi_x\cap T_x\Sigma, x\in \partial \Sigma,$ along $\partial \Sigma$ that points out
of $\Sigma,$ then there is an Euler class of $\xi|_\Sigma$ relative to $v.$ 
\begin{thm}\label{thm:bounds}
Let $M$ be a closed oriented irreducible 3--manifold and $\xi$ a plane field. Let $\Sigma$ be a surface embedded
in $M.$
\begin{enumerate}
\item (Thurston 1986, \cite{Thurston86}) If $\xi$ is a Reebless foliation and $\Sigma$ is closed, then 
\[
|\langle e(\xi),[\Sigma]\rangle|\leq -\chi(\Sigma),\quad \text{if } \Sigma\not=S^2
\]
and otherwise 
\[
|\langle e(\xi),[\Sigma]\rangle|=0.
\]
If $\Sigma$ has boundary transverse to $\xi,$ then
\[
\langle e(\xi), [\Sigma]\rangle \leq -\chi(\Sigma).
\]
\item (Eliashberg 1992, \cite{Eliashberg92a}) If $\xi$ is a tight positive contact structure and $\Sigma$ is closed, then 
\[
|\langle e(\xi),[\Sigma]\rangle|\leq -\chi(\Sigma),\quad \text{if } \Sigma\not=S^2
\]
and otherwise 
\[
|\langle e(\xi),[\Sigma]\rangle|=0.
\]
If $\Sigma$ has boundary transverse to $\xi,$ then
\[
\langle e(\xi), [\Sigma]\rangle \leq -\chi(\Sigma).
\]
\end{enumerate}
\end{thm}

\bhw
Show that this theorem implies there are only finitely many elements in $H^2(M;\Z)$ that can be the
Euler class of a Reebless foliation or tight contact structure. 
\ehw

The easiest way to prove a contact structure is tight is to find a symplectic filling. Recall a 4--manifold $X$ is a
\dfn{symplectic manifold} if there is a 2-form $\omega$ such that $d\omega=0$ and $\omega\wedge \omega$ is 
a never zero 4-form. Note $\omega\wedge \omega$ is a volume form on $X.$ We always assume $X$ is oriented 
by this form. If $M=\partial X$ (as oriented manifolds) 
and $\xi$ is a positive contact structure on $M$ then we say that $\omega$ \dfn{dominates} $\xi$ if $\omega|_\xi>0,$
(by this we mean $\omega(v,w)>0$ for any oriented basis $v,w$ for $\xi$). 
\bbr
It is very important that $M=\partial X$ as oriented manifolds. However, if you are not too interested in the orientation on $\xi$ then 
the condition that $\omega|_\xi>0$ can be easily arranged 
if $\omega|_\xi\not=0$ by reversing orientation on $\xi$. 
\eer
If $(M,\xi)$ is one component of a contact manifold $(M',\xi')$ and $(X,\omega)$ is a compact symplectic manifold
for which $\omega$ dominates $\xi'$ then we say that $(X,\omega)$ is a \dfn{weak symplectic semi-filling} of
$(M,\xi).$ If $M'$ is connected (that is $M'=M$) then we say $(X,\omega)$ is a \dfn{weak symplectic filling} of $(M,\xi).$
The reason we have brought up weak symplectic fillings is the following theorem.
\begin{thm}[Gromov-Eliashberg, \cite{Eliashberg90b, Gromov85}]
If $(M,\xi)$ is a weakly symplectically semi-fillable contact structure then $\xi$ is tight.
\end{thm}
\bex
Consider $S^3$ as the unit sphere in $\C^2=\R^4.$ On $\C^2$ we use polar coordinates $(r_1,\theta_1,r_2,\theta_2).$
Let $\alpha=r_1^2\, d\theta_1 + r_2^2\, d\theta_2.$ Earlier you checked that $\alpha$ restricted to $S^3$ is
a contact form for $\xi=\ker \alpha.$ Clearly,
\[
\omega = d\alpha= 2r_1\, dr_1\wedge d\theta_1 + 2r_2\, dr_2\wedge d\theta_2
\]
is a symplectic form on $\C^2.$ So in particular, it is a symplectic form on $B^4.$ Note
\[
\omega|_\xi=d\alpha|_\xi>0
\]
(since $\alpha\wedge d\alpha>0$).  Thus $(B^4,\omega)$ is a weak symplectic filling of $(S^3,\xi),$ and hence
$\xi$ is tight.
\eex
\bbr
This result should be compared to the famous result of Novikov that any foliation on $S^3$ must have Reeb components, \cite{Novikov65}.
(Recall, tight contact structures have no Lutz tubes and no overtwisted disks.)
\eer
We will find many other fillable, and hence tight,  contact structures later, but now the obvious question is: Are all
tight contact structures fillable? The answer is No.
\begin{thm}[Etnyre-Honda 2002, \cite{EtnyreHonda02b}]
There are tight but not weakly semi-fillable contact structures.
\end{thm}
This ends our brief introduction to contact structures and foliations. To learn more about foliations see \cite{CandelConlon00, CandelConlon03, HectorHirsch87} and to learn more about contact structures see
\cite{McDuffSalamon98}. We now move on to the perturbations of 
foliations into contact structures. 

\section{From foliations to contact structures}\label{sec:f2c}
Consider the interesting foliation on $S^2\times S^1$ given by 
\[
\zeta=\ker \alpha,
\]
where $\alpha=d\theta$ and $\theta$ is the coordinate on the $S^1$ factor. So $\zeta$ is given by the tangents to the
two spheres $S^2\times\{\theta\}.$
\begin{thm}[Eliashberg and Thurston, 1998 \cite{EliashbergThurston98}]\label{thm:f2c}
Any oriented $C^2$-foliation $\xi$ on an oriented 3-manifold $M,$ other than the foliation $\zeta$ of $S^2\times S^1,$
may be $C^0$-approximated by a positive and a negative contact structure. 
\end{thm}
To make sense of this theorem we need to have a topology on the space of plane fields. To this end recall that we can associate to the tangent bundle of $M$ the bundle of 2-planes. In other words at each point of $p\in M$ we replace $T_pM=\R^3$ with the Grassmann of 2-planes in $\R^3,$ which we denote $G_{2,3}.$
Sections of this new bundle are equivalent to plane fields on $M.$ Thus when we say we have  a $C^k$-plane field that means the corresponding section is 
$C^k$-smooth. Moreover, on the space of sections we have the topology of $C^k$-convergence, so we can talk about $C^k$-neighborhoods of plane fields. 
\bhw
Reinterpret the space of sections and the topology on them in terms of 1-forms. 
\ehw

We say a foliation $\xi$ can be \dfn{$\C^k$-deformed} into a contact structure if there is a $C^k$-family $\xi_t, t\in[0,\epsilon],$
so that $\xi_0=\xi$ and $\xi_t$ is a contact structure for $t>0.$ We also say $\xi$ can be \dfn{$C^k$-approximated} by a 
contact structure if in any $C^k$-neighborhood of $\xi$ there is a contact structure. It is clear that if $\xi$ can be $C^k$-deformed
into a contact structure then it can also be $C^k$-approximates by one too. 
\bex
Consider $T^3$ thought of as $\R^3,$ with coordinates $(x,y,z),$ modulo the action of the lattice $\Z^3.$ Now consider the 1-form
\[
\alpha_n^t= dz + t ((\cos\, 2\pi z)\, dx + (\sin\, 2\pi z)\, dy),
\]
where $n$ is any positive integer.
When $t=0$ we get $\alpha_n^0=dz$ which defines the foliation of $T^3$ by constant $z, $ $T^2$'s. When $t>0$ we get positive
contact structures $\xi_n^t.$ Note that by Gray's theorem if we fix $n$ then all the $\xi_n^t$ are isotopic for $t>0.$ So we
can unambiguously talk about $\xi_n$ (no $t$ dependence).  It is a result of Kanda \cite{Kanda97} and Giroux \cite{Giroux00} that the $\xi_n$ are
all distinct and, up to contactomorphism, give all tight contact structures on $T^3.$ When $t<0$ then note the $\xi_n^t$ are negative contact
structures on $T^3.$ 
\eex
\bbr
We make a few observations about the theorem and this example.
\begin{enumerate}
\item This last example shows that a fixed foliation can be approximated by (and even deformed into) infinitely many different 
contact structures! (We are not claiming that any foliation can be approximated by infinitely many different contact structures, just that 
some can.) This is somewhat surprising, as you might think that if contact structures are sufficiently close then you can deform
one to the other through contact structures and hence they would be isotopic/contactomorphic. This example demonstrates that this is not
the case.
\item The theorem only gives a contact approximation to the foliation not a deformation. It is possible the theorem is true with ``approximation''
replaced by ``deformation''. 
\item We lose smoothness in the theorem. We must start with a $C^2$-foliation, but the approximation is only $C^0$-close to the original foliation. 
It is possible that the theorem is true with ``$C^0$'' replaced with ``$C^2$''.
\end{enumerate}
\eer

Before we begin to sketch the proof of Theorem~\ref{thm:f2c} we first consider why the foliation $\zeta$ on $S^2\times S^1$ is so special.
This is indicated in the following two theorems. 
\begin{thm}[Reeb stability for confoliations, Eliashberg and Thurston, 1998 \cite{EliashbergThurston98}]
Suppose a confoliation $\xi$ on $M$ admits an embedded integral 2-sphere $S$ ({\em i.e.} for all $x\in S$ we have $T_xS=\xi_x$) then
$(M,\xi)$ is diffeomorphic to $(S^2\times S^1, \zeta).$
\end{thm}
This theorem is well known, and easier to prove, when $\xi$ is a foliation.
\bhw
Try to prove this theorem under the assumption that $\xi$ is a foliation.\\ Hint: Try to show the subset of $M$ that is foliated by $S^2$'s is
both open and closed. If you are having trouble maybe read ahead a few pages and come back and try again. 
\ehw
\begin{thm}[Eliashberg and Thurston, 1998 \cite{EliashbergThurston98}]
There is a $C^0$-neighborhood of $\zeta$ such that any confoliation of $S^2\times S^1$ in that neighborhood is 
diffeomorphic to $\zeta.$
\end{thm}
This last theorem really explains why we must have the exceptional case of $(S^2\times S^1, \zeta)$ in Theorem~\ref{thm:f2c}. We can
also use the Reeb stability theorem to see how not to try to perturb a foliation into a contact structure. Indeed, the next theorem 
implies that you can not ``locally'' perturb a foliation into a contact structure.
\begin{thm}
Let $\xi$ be a confoliation on the 3-ball $B$ which is standard near $\partial B$ ({\em i.e.} near $\partial B,$ $\xi$ is given by 
$\ker dz$). Then $\xi$ is a foliation and is diffeomorphic to the standard foliation on $B.$
\end{thm}

\begin{proof}
Given $(B,\xi)$ we can find an embedding of $B$ into $S^2\times S^1$ so that near $\partial B$ $\xi$ agrees with $\zeta.$ Now define
$\zeta'$ to be $\zeta$ on $(S^2\times S^1)\setminus B$ and $\xi$ on $B.$ This is a confoliation on $S^2\times S^1$ that has an integral 
sphere (we can certainly choose our initial embedding of $B$ so that it misses $S^2\times \{pt\}$ for some $pt$). Thus the Reeb stability 
theorem implies that $\zeta'$ is diffeomorphic to $\zeta.$ This diffeomorphism shows that $\xi$ is a foliation diffeomorphic to
$\zeta$ restricted to the image of the ball.
\bhw
Show that $\xi$ is indeed standard the standard foliation on $B^3.$ 
\ehw
\end{proof}

We can break the proof the main theorem, Theorem~\ref{thm:f2c}, into two parts:
\begin{enumerate}
\item[Part 1] Perturb $\xi$ into a confoliation $\xi'$ such that $\xi'$ is contact on a ``sufficiently large'' part of $M.$
\item[Part 2] Perturb $\xi'$ into a contact structure. 
\end{enumerate}
These steps seem overly simplistic, but this is a good outline of the strategy. We will actually start with Part 2 so that we can figure out what 
``sufficiently large'' in Part 1 actually means! During the proof of Theorem~\ref{thm:f2c} we will find a positive contact structure approximating
$\xi.$ The proof for a negative contact structure is similar. From this point on when we say ``contact structure'' we mean ``positive contact structure''.

\subsection{Part 2 of the proof of Theorem~\ref{thm:f2c}}
Suppose we are given a confoliation $\xi'$ on a closed oriented manifold $M,$  then set
\[
H(\xi')=\{x\in M :  \xi'_x \text{ is contact at } x, \text{ ({\em i.e.} }(\alpha\wedge d\alpha)_x>0 \text{)}\}.
\]
The set $H(\xi')$ is called the \dfn{hot zone} or \dfn{contact region}. (The reason for the terminology ``hot zone'' will be clear from the
discussion below.) Now set
\begin{align*}
G(\xi')=\{ x\in M : &\text{ there is a path } \gamma \text{ from $x$ to $y$} \\& \text{ such that }  y\in H(\xi') \text{ and } \gamma \text{ is tangent to } \xi'\}.
\end{align*}

The precise statement of Part 2 of the proof of Theorem~\ref{thm:f2c} is contained in the following theorem. 
\begin{thm}\label{step2}
If $G(\xi')=M$ then $\xi'$ can be $C^\infty$-deformed into a contact structure.
\end{thm}
There are two proofs of this theorem. An {\em analytic proof} (due to Altschuler \cite{Altschuler95}) and a {\em topological proof} (due to Eliashberg
and Thurston \cite{EliashbergThurston98}). We sketch both these approaches.

\noindent 
{\bf The analytic way:} Choose a Riemannian metric on $M$ and further choose a 1-from $\alpha$ such that $\xi' =\ker \alpha$ and
$|\alpha|=1$ at all points of $M.$ Now consider the equations
\begin{align*}
\frac{\partial}{\partial t} \beta &= *(\alpha\wedge df),\\
\beta_0&=\alpha,
\end{align*}
where $f=*(\alpha\wedge d\beta +\beta\wedge d\alpha)$ and $\beta$ is a section of $T^*M\times \R^+.$ Here $\R^+$ is the non-negative
real numbers and we think of $\beta$ as a 1-parameter family of 1-forms on $M.$ (We will frequently denote the time dependence of $\beta$
as a subscript.) In these equations we are given $\alpha$ and we
are trying to solve for $\beta.$

These equations are a weakly-parabolic system and Altschuler proved that given $\alpha$ as above there is a unique smooth solution 
for $t\in [0,\infty).$ See \cite{Altschuler95}. 

The function $f$ also evolves by a weakly-parabolic equation:
\[
\frac{\partial f}{\partial t} = \Delta_\alpha f + \nabla_X f,
\]
where $X$ is some time dependent vector field and $\Delta_\alpha$ is the ``Laplacian on $\ker \alpha$''. Intuitively $\Delta_\alpha$ is the sum of partial 
derivatives in the directions tangent to $\xi.$ For the precise definition see \cite{Altschuler95}.
A version of the maximum principle for this equation implies the following fundamental property:
 If $q$ is connected to a point $p$ by a path tangent to $\xi'$ and $f(p,0)>0,$ then $f(q,t)>0$ for all $t>0.$ In particular,
 if $f(p,0)>0$ then $f(p,t)>0$ for all $t>0.$

As the archetypical parabolic equation is the heat equation we might think of $f$ as representing heat then this last property says that
``Heat flows infinitely fast to all points of $M$ accessible to the hot zone''. Note that for any point $p\in H(\xi')$ we have $f(p)>0.$ 
So if $G(\xi')=M$ then $f(p,t)>0$ for all $t>0.$ This is good
because if we set 
\[
\eta = \alpha + \epsilon \beta_1
\]
then
\[
d\eta = d\alpha + \epsilon d\beta_1,
\]
and 
\[
\eta \wedge d\eta = \alpha\wedge d\alpha + \epsilon (\alpha\wedge d\beta_1 + \beta_1\wedge d\alpha) + \epsilon^2 \beta_1\wedge d\beta_1.
\]
Since $\xi'$ is a confoliation we know the first term on the right hand side is greater than or equal to 0. The second term is just $\epsilon (*f_1)$
which is positive everywhere since $G(\xi')=M.$ Finally, the last term might have any sign, but by choosing $\epsilon$ small enough its 
magnitude will be smaller than that of the second term. Thus $\eta$ is a contact form for all $\epsilon$ near zero and we have constructed 
the deformation of $\xi'$ from Theorem~\ref{step2}.

\noindent
{\bf The topological way:}  This proof is due to Eliashberg and Thurston, \cite{EliashbergThurston98}.
Given an arc $\gamma$ tangent to $\xi'$ with one end in a contact region, we want to show how to extend the contact region to contain a neighborhood
of $\gamma.$ To this end parameterize $\gamma$ by $[0,1]$ and we find a neighborhood $N$ of $\gamma$ of the form $N=[-1,1]\times [0,1]\times [-1,1],$ with
coordinates $(x,y,z),$ such that 
\begin{enumerate}
\item $\gamma=\{x=0,z=0\}$,
\item $\xi'=\ker \alpha$ where 
\[
\alpha=dz-a(x,y,z)\, dx, \text{ and}
\]
\item  $\xi'$ is contact near $\{y=1\}.$
\end{enumerate}
(A slight modification of your argument in the proof of Lemma~\ref{lem:locform} above will suffice to find this neighborhood.)
The main goal now is to prove the following lemma.
\begin{lem}
There is a $C^\infty$-deformation of $\xi'$ supported in $N$ to a plane field that is  contact on the interior of $N.$
\end{lem}
\begin{proof}
Since $\xi'$ is a confoliation we know, from Lemma~\ref{lem:type}, that $\frac{\partial a}{\partial y}\geq 0$ in $N$ and near $\{y=1\}, \frac{\partial a}{\partial y}>0.$
Choosing $x_0$ and $z_0$ then  $a(x_0,y,z_0)$ is a non-decreasing function that is strictly increasing near $y=1.$ 
It is easy to choose
a new function of $y$ that is strictly increasing for all $y$ in $(0,1)$ and arbitrarily close to $a(x_0,y,z_0).$ The trick now is to simultaneously  perturb $a(x_0,y,z_0)$
for all $x_0$ and $z_0.$ To this end, choose $\delta$ so that all the $a(x_0,y,z_0)$ are strictly increasing for $y$ in $[1-\delta, 1].$ Now let $m$ be the minimum of
$\frac{\partial a}{\partial y}(x_0,y,z_0)$ for all $(x_0,y,z_0)\in [-1,1]\times [1-\delta,1]\times [-1,1].$ Choose a positive 
function $f(y)$ so that $f(y)$ is strictly increasing on 
$[0,1-\delta]$ and decreasing on $[1-\delta,1],$ but with derivative larger than $-m.$ Moreover, we can choose $f$ so that $f(y)$ is bounded by a constant that is a small as we like. Choose a cut-off function $g(w)$ such that $g(w)$ is 1 for $w\in[0,1-\epsilon],$
for $\epsilon$ small, and $g(w)$ and all its derivatives are zero at 1. If we 
set
\[
\tilde{a}_t(x,y,z)=a(x,y,z) + tg(x^2)g(z^2) f(y),
\]
then one may easily check that $\tilde{\alpha}_t=dz-\tilde{a}_t(x,y,z)\, dx$ is the desired deformation of $\alpha.$
\end{proof}
We now know we can turn $\xi'$ into a contact structure in neighborhoods of arcs tangent to $\xi'.$ We can find a finite number $\gamma_1,\ldots, \gamma_n$
of arcs tangent to $\xi'$ that start in the contact region and so that each has a neighborhood $N_i$ as above and the $N_i$'s cover $M\setminus H(\xi').$ We would now like
to change $\xi'$ into a contact structure in each of the $N_i$'s, but notice that when we do this for $N_1$ we have changed $\xi'$ on some of the other $N_i$'s
so the form of the 1-form representing $\xi'$ on these other neighborhoods changes. So we might not be able to continue on the other $N_i$'s.
\bhw
Show that  if the perturbation on $N_1$ is sufficiently small then we can slightly modify the $N_i$'s so that we still have the appropriate form to apply the lemma.
\ehw

\subsection{Part 1 of the proof of Theorem~\ref{thm:f2c}}

Now we know that given a foliation $\xi$ if we can perturb it to through confoliations to $\xi'$ so that $G(\xi')=M,$ that is every point in $M$ can be connected to a 
a ``contact region'' by a path tangent to $\xi'$, then we can move on to Part 2 and perturb $\xi'$ into a contact structure. So we now want to figure out how to do this
first perturbation. Specifically we want to see how to prove the following precise statement of Part (1) of the proof of Theorem~\ref{thm:f2c}.
\begin{lem}\label{lem:part1}
Any $C^2$-foliation $\xi$ on an oriented 3-manifold $M,$ other than the foliation $\zeta$ on $S^2\times S^1,$ can be $C^0$-approximated by a
confoliation $\xi'$ such that $G(\xi')=M.$
\end{lem}
Note that in Part 2 of the proof we $C^\infty$-deformed $\xi'$ into a contact structure. So it is in this part of the proof that we will see the loss of
smoothness and the fact that we can only approximate $\xi$ by a contact structure and not necessarily deform $\xi$ into a contact structure. 

\begin{figure}[ht]
  \relabelbox \small {\centerline{\epsfbox{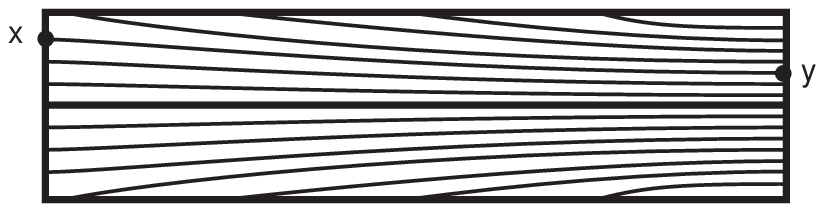}}} 
  \relabel{x}{$x$}
  \relabel{y}{$y=\phi(x)$}
  \endrelabelbox
        \caption{The foliation induced on the annulus $A$ and the map $\phi.$ The heavier line in the center is $\gamma.$}
        \label{fig:holo}
\end{figure}

The main tool we need to introduce contact regions into $\xi$ is {\em holonomy}. Fix a foliation $\xi$ on $M$ and a closed oriented curve $\gamma$ in $M$ that is tangent
to $\xi.$ We can embed an annulus $A=(-\epsilon,\epsilon) \times S^1$ in $M$ such that 
\begin{enumerate}
\item $\{0\}\times S^1=\gamma$ (as oriented curves),
\item $A$ is transverse to $\xi,$ and
\item $(-\epsilon, \epsilon)\times\{p\}$ is transverse to $\xi$ for all $p\in S^1.$
\end{enumerate}
Thus $\xi$ induces a foliation on $A,$ which we denote $A_\xi.$ Clearly $\gamma$ is a leaf in this foliation. 
If we restrict the projection map $\pi:A\to S^1$ to a leaf $l$ of $A_\xi$ then $\pi$ will be a local diffeomorphism. Thus if we choose a point $x$ on a $l$ and let $C$ be a 
path in $S^1$ based at $\pi(x)$ then we may try to lift $C$ to a path in $l.$ Just like when you prove path lifting for covering spaces, it is easy to see you can lift $C$
as long as the lift does not  ``run off of $A$''. For example if $C$ just parameterizes $\gamma,$ and we lift $C$ starting at a point $x\in \{0\}\times S^1$ then 
the lifted curve $C$ just parameterizes $\{0\}\times S^1.$ Similarly if $x$ is a point very near $\{0\}\times S^1$ then $C$ will lift to a closed leaf of $A_\xi$ or
an arc contained in a non-colsed leaf. If we fix a point $p$ in $S^1$ and set $I=(-\epsilon, \epsilon)\times \{p\}$ then we have just seen for $x$ near $0$ in $I,$
we can lift $C$ to a curve $C_x$ tangent to a leaf of $A_\xi.$ If $x$ is close enough to $0$ then 
the curve $C_x$ will intersect $I$ in either one point $\{x\}$ (if $C_x$ parameterizes a closed leaf) or
two points $\{x, y\}.$ 
See Figure~\ref{fig:holo}.
So we can define a map
\[
\phi_\gamma\co I'\to I
\]
where $I'$ is a sub-interval of $I$ containing 0. This map sends $x\in I'$ to $x$ if $C_x$ is a closed curve and $y$ otherwise. Note $\phi_\gamma(0)=0.$
This map is called the \dfn{holonomy along $\gamma$}. Actually, to precisely define the holonomy along $\gamma$ takes more work, in particular we need to make
$\phi_\gamma$ independent of $A.$ While it is not so important for us here, let's consider some of the subtleties. First, if we change $A$ a little then the map
$\phi_\gamma$ may be conjugated by a diffeomorphism of $I.$ Secondly, when we change $A$ the interval $I'$ on which the map $\phi_\gamma$ is well-defined
might change. So it is really the conjugacy class of the germ of $\phi_\gamma$ that is the holonomy of $\xi$ along $\gamma.$ While these subtleties are important
we will largely ignore them. The careful reader should make sure we don't miss anything by doing this. Also note that given a leaf $l$ of $\xi$ we get a map
\[
\Phi\co \pi_1(l)\to \{\text{germs of maps of } I \text{ to itself at } 0\},
\]
by sending $\gamma\in\pi_1(l)$ to $\phi_\gamma.$
The (conjugacy class) of the image of $\Phi$ is called the holonomy group of $l.$ 

We will call the holonomy along a curve $\gamma,$ \dfn{non-tivial} if $\phi_\gamma\not = id_I$ and \dfn{linearly non-trivial} if $\phi'_\gamma(0)\not = 1.$ The 
holonomy will be said to be \dfn{attracting} (\dfn{repelling}) if $|\phi_\gamma(x)|< |x|$ ($|\phi_\gamma(x)|>|x|$) for all $x,$ near 0. Finally, we say the 
holonomy is \dfn{sometimes attracting (repelling)} if $|\phi_\gamma(x)|< |x|$ ($|\phi_\gamma(x)|>|x|$) for a sequence of $x$ approaching 0 from both sides
(this is the same as saying for $x$ on intervals arbitrarily close to, but not necessarily including, 0).

So why are we interested in holonomy? It helps us create contact regions!
\begin{thm}\label{holgood}
Let $(M,\xi)$ be a $C^k$-foliated manifold. 
\begin{enumerate}
\item If $\gamma$ is a curve tangent to $\xi$ and has non-tirivial linear holonomy, then there are neighborhoods $N$ and $N'$ of $\gamma$
such that $N\subset N'$ and $\xi$ can be $C^k$-deformed through confoliations so that it is a positive contact structure in $N,$ unchanged outside of $N'$ and
diffeomorphic to $\xi$ outside $N.$
\item If $\gamma$ is a curve tangent to $\xi$ and has sometimes attracting (repelling) holonomy, then there are neighborhoods $N$ and $N'$ of $\gamma$
such that $N\subset N'$ and $\xi$ can be $C^0$-approximated by a confoliation so that it is a positive contact structure in $N,$ unchanged outside of $N'$ and
diffeomorphic to $\xi$ outside $N.$
\end{enumerate}
\end{thm}

\begin{proof}
Part (1) of the theorem can be proved in two ways. We illustrate both ways, since the first is simpler and the second generalizes to provide a proof of part (2).
For the first proof let $U=\gamma\times[-1,1]\times[-1,1]$ be a neighborhood of $\gamma$ with coordinates $(x,y,z)$ corresponding to the three factors in the
product. From Lemmas~\ref{lem:locform} and \ref{lem:type} these coordinates can be chosen so that there is a 1-form $\alpha$ satisfying 
$\xi=\ker\alpha$ and
\[\alpha =  dz - a(x,z)\, dx.\]
We in addition claim that, because $\gamma$ has linear holonomy, the coordinates can be chosen so that 
\[
-\frac{\partial a}{\partial z}\geq C,
\]
for some positive constant $C>0.$
\bhw
Justify this claim.\\
Hint: Construct a model situation satisfying the inequality and then make sure this can be 
inserted into the given situation. More specifically,  on $\gamma\times[-1,1]$ (where we think of $[-1,1]$ as the $z$-direction) try to 
abstractly construct an $a(z)$ with the desired property and inducing the
same holonomy map as in the given situation. Then use the ideas in the proof of Lemma~\ref{lem:locform} to construct the neighborhood $U$ above.
\ehw

Now let $h\colon [0,1]\to \R$ be a function so that $h(0)=1, h(1)=0$ and $h$ is decreasing. Set
\[
\beta=h(y^2+z^2)\, dy.
\]
Note that
\[
\alpha\wedge d\beta + \beta\wedge d\alpha= (-a_zh+2zah')\, dx\wedge dy\wedge dz.
\]
Since $a$ is strictly 
decreasing with respect to $z$ and equal to $0$ at $z=0$ we know that $(-a_zh+2zah')>0.$ Thus if $\alpha_\epsilon=\alpha+\epsilon \beta$
then one may easily check that $\alpha_\epsilon$ is a contact form on $N=\gamma\times \{(y,z)|y^2+z^2< 1\}$ for $\epsilon>0$ and $\alpha_\epsilon =\alpha$
outside of $N.$ Thus we have completed the proof of part (1) of the theorem.

Now for the second proof of (1), that will generalize to a proof of (2). Again consider a neighborhood $U=\gamma\times[-1,1]\times[-1,1]$ as above. Let 
\[
A_{y_0}=\{(x,y,z)\in U| y=y_0\}.
\]
Note each annulus $A_{y_0}$ 
has a foliation induced on it by $\xi.$ Moreover $\xi$ at any point in $U$ is spanned by a vector in the foliation on $A_{y_0}$ and $\frac{\partial}
{\partial y}.$ Thus the foliations induced on the $A_{y_0}$'s essentially determine $\xi.$ We will describe a way to change the foliations on the annuli $A_{y_0}$
to get a contact structure. 

The idea will be to construct a diffeomorphism $F\colon U\to U$ such that 
\begin{enumerate}
\item $F$ is the identity on $\partial U$ (and $C^k$-tangent to the identity on $\partial U$),
\item $F$ preserves the annuli $A_{y_0}$ and 
\item the slope of the foliation on $A_{y_0},$ for $y_0\in (-1,1),$
after pushing the foliation forward by the diffeomorphism is always less than the slope before the diffeomorphism. 
\end{enumerate}
See Figure~\ref{fig:twist}. 
\begin{figure}[ht]
  \relabelbox \small {\centerline{\epsfbox{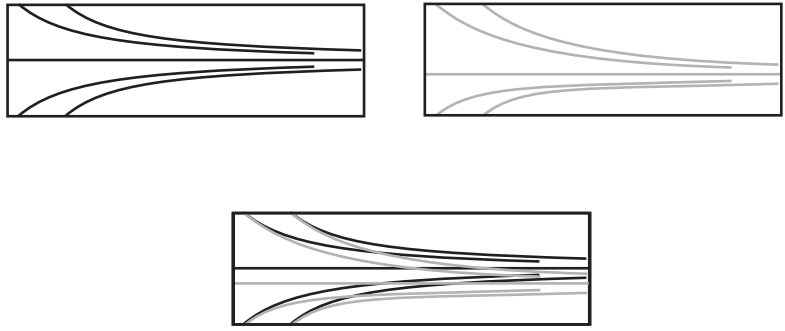}}} 
  \endrelabelbox
        \caption{The foliation on $A_y$ induced by $\xi$ on the left and the foliation on $A_y$ after applying the diffeomorphisms on the right. On the
        	bottom are the two foliations superimposed on each other. Here one can see that the slope of the foliation at any point (away form the boundary)
	 after the diffeomorphism is less than the slope of the original foliation.}
        \label{fig:twist}
\end{figure}
Now we have two foliations on $U.$ The first is the one we started with $\xi$ and the second is $\xi'=F_*(\xi).$ 
Note that both $\xi$ and $\xi'$ contain $\frac{\partial}{\partial y}$ so they are determined by the foliations induced on the $A_{y}$'s. Also note that $\xi=\xi'$ 
on $\partial U.$ Thus we could replace $\xi|_U$ by $\xi'$ and we would have another foliation on $M.$ Instead we construct $\widetilde{\xi}$ on $U.$ Let $U_{y_0}$
denote the set $\{(x,y,z)| y\leq y_0\}.$ We define $\widetilde{\xi}$ to be equal to $\xi$ on $U_{-\frac 12}$ and equal to $\xi'$ on $U\setminus U_{\frac 12}.$ On the 
remaining region we will interpolate between these foliations. To define $\widetilde{\xi}$ in this region fix a point $(x_0,z_0)$ and consider the line
segment $l_{(x_0,z_0)}=\{(x,y,z)|x=x_0, y\in [-\frac 12, \frac 12], z=z_0\}.$ The slope $s_{-\frac 12}$
of $\widetilde{\xi}$ on $A_{-\frac 12}$ at $(x_0,z_0)$ is greater  than the slope $s_{\frac 12}$ of
$\widetilde{\xi}$ on $A_{\frac 12}$ at $(x_0,z_0).$ Thus we can define $\widetilde{\xi}$ to be the plane field that is always tangent to $l_{(x_0,z_0)}$ and the slope
of the intersection with $\widetilde{\xi}$ with $A_{y}$ uniformly decreases from $s_{-\frac 12}$ to $s_{\frac 12}$ as $y$ goes from $-\frac 12$ to $\frac 12.$
From Lemma~\ref{lem:type} it should be clear that $\widetilde{\xi}$ is a positive contact structure in $U_{\frac 12}\setminus U_{-\frac 12}.$
\bhw
Write down forms for $\xi$ and $\xi'.$ Using the properties of $F$ listed above write down a form for $\widetilde{\xi}$ and show that this is a contact form 
in the region  $U_{\frac 12}\setminus U_{-\frac 12}.$
\ehw
So we are left to construct the diffeomorphism $F.$ To this end we need a lemma.
\begin{lem}
Let $v_x$ be a family of smooth functions on $[-1,1]$ such that $v_x(0)=0$ and $v_x$ are monotonically increasing for all $x.$ Then there exists a diffeomorphism
$f\colon [-1,1]\to [-1,1],$ $C^\infty$-close to the identity and $C^\infty$-tangent to the identity at $\{-1,1\}$ and satisfying 
\[
f'(z) v_x(z)> v_x(f(z)),
\]
for all $z\in (-1,1)$ and $x.$
\end{lem}
\bhw
Prove this lemma.
\ehw
Now in the case when $\gamma$ has linear holonomy we can choose our coordinates as in the first proof of this theorem above. Thus if we set
$v_x(z)=-a(x,z)$ then $v_x$ satisfies the conditions of the lemma and we thus get the function $f$ from the lemma. Let $f_s, s\in [-1, -\frac 12]$ be a linear homotopy 
between $f_{-1}$ the identity on $[-1,1]$ and $f_{-\frac 12}=f.$ If we choose this isotopy to be tangent to the identity at the end points we can extend it to be the 
constant isotopy for $s\in[-\frac 12, 1].$ Finally we can set  $F(x,y,z)=(x,y, f_y(z)).$
\bhw
Check that $F$ has the properties listed above.
\ehw
Thus we have finished our second proof of the theorem in the case of linear holonomy along $\gamma.$ For the case of sometimes attracting (repelling) holonomy
we just need to see how to construct the diffeomorphism $F.$ For this we need the following lemma.
\begin{lem}
Let $v_x$ be a smooth family of $C^1$-functions on $(-1,1)$ such that $v_x(0)=0$ and there is a sequence of points $z_n$ decreasing to $0$ and points $z'_n$
increasing to $0$ such that $v_x(z_n)>0$ and $v_x(z'_n)<0$ for all $x$ and $n.$ Then for any sufficiently small $\epsilon>0$ there exists a diffeomorphism 
$f\colon (-1,1)\to(-1,1)$  which is fixed outside the interval $(-\epsilon, \epsilon)$ and satisfies 
\[f'(z)v_x(z)>v_x(f(z)),\]
for all $z\in(-\epsilon, \epsilon)$ and $x.$
\end{lem}
\bhw
Prove this lemma.
\ehw
\bhw
Find a suitable modification of the argument in the linear holonomy case to handle this more general case.
\ehw
\end{proof}

OK, so now we know that we can use holonomy to create regions of contact in a foliation. We are now ready to prove Lemma~\ref{lem:part1}.
\begin{proof}[Proof of Lemma~\ref{lem:part1}]
We need perturb $\xi$ to a confoliation $\xi'$ so that any point in $M$ can be connected to a region where $\xi'$ is contact by a path tangent to $\xi'.$ From Theorem~\ref{holgood} this is the same as proving that we can perturb $\xi$ into another foliation so that every leaf in the new foliation is arbitrarily close to a
leaf with a curve having holonomy. To this end we study minimal sets. A \dfn{minimal set} in a foliation $\xi$ is a non-empty closed union of leaves that contains 
no smaller such set. 
\bhw
Show that a set is minimal if it is a non-empty closed union of leaves and every leaf in the set is dense in the set.
\ehw
\bhw
Show that every leaf in $\xi$ limits to some minimal set. 
\ehw
Thus we just need to see that we can perturb $\xi$ so that every minimal set has a curve with holonomy. Minimal sets can be quite complicated for arbitrary 
foliations, but if we restrict to $C^2$ foliations or better we have the following classification of minimal sets.
\begin{thm}
In a $C^2$-foliation every minimal set is either
\begin{enumerate}
\item all of $M$ (in which case the foliation is called \dfn{minimal}),
\item a closed compact leaf, or
\item an \dfn{exceptional minimal set}.
\end{enumerate}
\end{thm}
For more discussion of this theorem and exceptional minimal sets see \cite{HectorHirsch87}.
For our purposes an exceptional minimal set will mean a minimal set that is not of type (1) or (2). With this definition the above theorem is not too
hard to prove! However, there is a lot of structure to exceptional minimal sets. We will need
\begin{thm}[Sachsteder, 1965 \cite{Sacksteder65}]
Exceptional minimal sets contain leaves with linear holonomy.
\end{thm}
Thus from our perspective all we need to know is that we can always find a curve with linear holonomy in an exceptional minimal set. We have
just two other types of minimal sets to consider. If $\xi$ is minimal, that is $M$ is the only minimal set of $\xi$ then there are two cases to consider. 
The first is when there is some holonomy in $\xi.$ In this case a result of Ghys (see \cite{EliashbergThurston98}) 
says that $\xi$ has linear holonomy as well. Thus we are done in this case.
If $\xi$ has no holonomy then a theorem of Tischler \cite{Tischler70} imples that $\xi$ my be $C^0$ approximated by a fibration over the circle. That is $M$ is a fibration
over the circle and $\xi$ is $C^0$-close to the foliation of $M$ by the fibers of this fibration. Note, by assumption, the fiber of this fibration cannot be $S^2.$ We claim
we can now $C^\infty$-approximate this foliation by surfaces by a $C^0$-foliation with only two closed leaves and linear holonomy along the closed leaves. 
We can then use the linear holonomy to create regions of contact and then finally perturb this confoliation into a contact structure, so the fact that at one point
the foliation was only $C^0$ is irrelevant (since the contact condition is an open condition on 1-forms, we can always perturb a contact 1-form to be as smooth as
we like). To perturb our fibration we can think of $M$ as a mapping torus
\[\Sigma\times[0,1]/\sim,\]
where $(0,x)\sim (1,f(x)),$ for some diffeomorphism $f$ of $\Sigma.$ We will perturb the foliation tangent to the fibers of this fibration so that $\Sigma\times\{0\}$ and
$\Sigma\times\{\frac 12\}$ are the only two closed leaves. To this end choose a separating curve $\gamma$ in $\Sigma$ and cut $\Sigma\times[0,\frac 12]$ open along
$\gamma\times[0,\frac 12].$ We can now shear all the leaves in the foliation as shown in Figure~\ref{fig:shear} and reglue. 
\begin{figure}[ht]
  \relabelbox \small {\centerline{\epsfbox{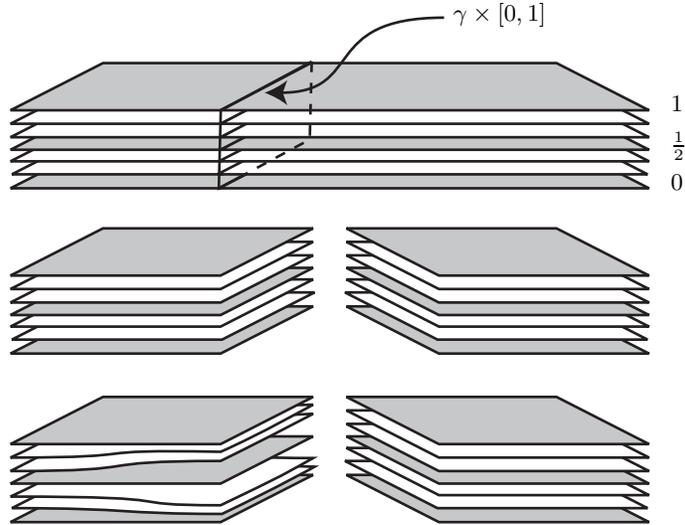}}} 
  \relabel{1}{$\frac 12$}
  \relabel{2}{$1$}
  \relabel{0}{$0$}
  \relabel{g}{$\gamma\times [0,1]$}
  \endrelabelbox
        \caption{The top figure is a neighborhood of $\gamma\times[0,1]$ in the foliation. In the middle figure the foliation has been cut open along $\gamma\times[0,1].$
        At the bottom we shear the leaves on one side of the cut open foliation. The gray leaves are $\Sigma\times \{0,\frac 12, 1\}.$}
        \label{fig:shear}
\end{figure}
If you would prefer a more rigorous description of this perturbation then consider $U=S^1\times[-1,1]\times [0,\frac 12]$ a neighborhood of $\gamma\times[0,\frac12].$ The
original foliation is by $S^1\times[-1,1]\times\{pt\}$'s. We can replace this foliation by $\ker\alpha$ where $\alpha=dz+g(y)h(z)dy$ and $h\colon [0,\frac 12]\to \R$ is
0 at $0$ and $\frac 12$ and negative on $(0,\frac 12)$ and $g\colon [-1,1]\to \R$ vanishes to high order at $-1$ and 1 and is positive on $(-1,1).$
Note we have not changed $\Sigma\times[0,\frac 12]$
but we have changed the foliation on it. In particular, if $\gamma'$ is a closed curve in $\Sigma\times\{\frac 12\}$ that intersects $\gamma$ in one point, then there will
be one-sided linear holonomy along $\gamma'.$
We can now do the same thing with $\gamma\times[\frac 12, 1]$ on $\Sigma\times[\frac 12,1].$ After this note $\gamma'$ in $\Sigma\times\{\frac 12\}$ now has
linear holonomy. If the monodromy $f$  of the surface bundle fixed $\gamma$ then we would have a nice smooth foliation with holonomy along $\gamma'$ in
$\Sigma\times\{0\}$ too. In general $\gamma$ will not be fixed by $f$ this leads to the lack of smoothness in our foliation, but it is not too hard to find a curve
in $\Sigma\times\{0\}$ that has linear holonomy. 
\bhw
Find this curve.\\
Hint: Find a $\gamma''$ for $f(\gamma)$ that acts as $\gamma'$ did for $\gamma.$ Then add  $\gamma''$ to $\gamma'.$
\ehw
We note that with a great deal more care this idea can be significantly generalized to produce a $C^2$-foliation by performing this shearing on 
infinitely many curves $\gamma.$ See \cite{EliashbergThurston98}. 

So we are done if $\xi$ is a minimal foliation. If $\xi$ is not minimal then we just need to worry about exceptional minimal sets and closed leaves. It turns out there
are always finitely many exceptional minimal sets, but this does not have to be the case for closed leaves. However we do have the following theorem.
\begin{thm}
The foliation $\xi$ can be $C^0$-perturbed so that it has only finitely many closed leaves.
\end{thm}
The idea to prove this theorem is very similar to the idea above. If you have a bunch of closed leaves you can try to shear as we did above. For all the details see
\cite{EliashbergThurston98}. From now on we assume there are only finitely many closed leaves. 

If our foliation has a closed leaf $\Sigma$ 
there are three cases to consider (1) no holonomy, (2) linear or weakly attracting/repelling holonomy and (3) holonomy but not weakly 
attracting/repelling. In case (1) we have the Reeb stability theorem \cite{HectorHirsch87} that says there is a neighborhood of the closed leaf foliated 
by closed leaves. This contradicts 
the fact that we arranged to have only finitely many closed leaves! Thus case (1) cannot happen. Case (2) is the good case where we can introduce regions of 
contact. In case (3) we know we have holonomy, 
but it is not weakly attracting/repelling.
If no curve on $\Sigma$ has holonomy on one side then a version of Reeb stability implies that $\Sigma$ has a one sided neighborhood that is foliated by closed
surfaces. Since we are assuming only finitely many surfaces this cannot happen.
Thus we can assume there is a curve $\gamma$ on the leave such that the holonomy is weakly attracting on one side
and weakly repelling on the other side. In this case we can split the manifold open along the leaf, glue in a product neighborhood $\Sigma\times [-\epsilon,\epsilon]$
and extend the foliation over this product just to be $\Sigma\times\{pt\}.$ Now use the shearing trick above to make sure that $\Sigma\times\{-\epsilon\}$ has 
weakly attracting holonomy along $\gamma$ and $\Sigma\times\{\epsilon\}$ has weakly repelling holonomy. See Figure~\ref{fig:insert}. 
\begin{figure}[ht]
  \relabelbox \small {\centerline{\epsfbox{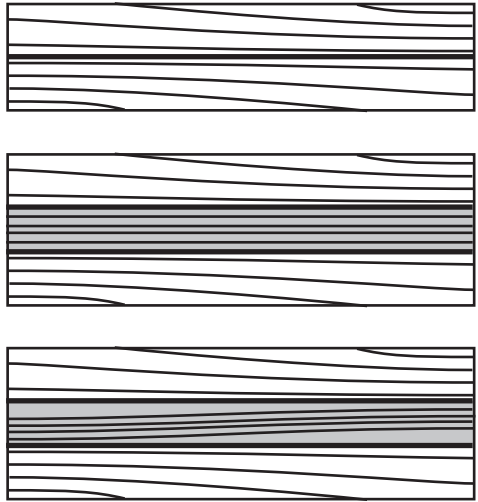}}} 
  \endrelabelbox
        \caption{The top figure is a transverse annulus about a curve $\gamma$ whose holonomy is attracting on one side and repelling on the other side. (The curve
        	$\gamma$ is the thicker line in the center.)
        	In the middle figure is the foliation on the transverse annulus after the foliation has been cut open and the trivial foliation on $\Sigma\times[-\epsilon, \epsilon]$
	is added. The bottom figure shows the final foliation after shearing has been done to create two copies of $\gamma$, one with attracting holonomy and the other
	with repelling holonomy.}
        \label{fig:insert}
\end{figure}

Putting this all together we see that after, possibly, perturbing our foliation we can assume there are finitely many exceptional minimal sets and finitely many closed
leaves and all these minimal sets contain leaves with linear or weakly attracting/repelling holonomy. Thus we can use Theorem~\ref{holgood} to perturb the foliation
into the desired confoliation.
\end{proof}

\section{Taut foliations and symplectic fillings}\label{taut2fill}
From Theorem~\ref{thm:f2c} we know that we can perturb a foliation $\xi$  into contact structures $\xi'$ 
(assuming of course that the foliation is not the trivial one on $S^1\times S^2$), but what can we say about $\xi'$? Is it tight? Is it fillable? Is it overtwisted?

To answer this recall that Thurston showed that if a foliation $\xi$ on $M$ is Reebless then for all surfaces $\Sigma\not= S^2$ embedded in $M$ we
have
\[
|\langle e(\xi),[\Sigma]\rangle|\leq -\chi (\Sigma).
\]
Thus if $\xi'$ is a contact structure $C^0$-approximating $\xi$ then $e(\xi')$ will satisfy the same inequality. This might lead one to believe that $\xi$ is tight,
but it does not constitute a proof that $\xi'$ is tight.
If we could arrange that $e(\xi')$ of the contact structure satisfied the inequality in Theorem~\ref{thm:bounds} involving transverse knots this would
be sufficient to conclude that $\xi'$ is tight. However, curves transverse to $\xi$ are different from curves transverse to $\xi'$ so we cannot conclude that $e(\xi')$
satisfies this inequality just because $e(\xi)$ does. None the less our intuition on this matter is indeed correct.
\begin{thm}
If $\xi'$ is a positive contact structure $C^0$-close to a Reebless foliation, then $\xi'$ is tight. Moreover, $\xi'$ is still tight when pulled back to the universal cover of
$M$.
\end{thm}
This result was originally stated in Eliashberg and Thurston's book \cite{EliashbergThurston98}, but an error in the proof was discovered in \cite{Colin02} which was corrected in \cite{EtnyreRLisT}.

It is interesting to note that there are many foliation that have Reeb components but still perturb to a tight contact structure. For example, 
recall that $S^3$ can be thought of as the union of two solid tori. We can foliate each of these tori with Reeb foliations. Topologically there are two ways to do this (depending on the direction of the spiraling of the leaves on both sides of the unique torus leaf). It turns out that one of these  
Reeb foliations on $S^3$ perturbs to an
overtwisted contact structure. The other Reeb foliation of $S^3$ perturbs to a tight contact structure. 
\bhw
Determine which foliation perturbs to an overtwisted contact structure. You might try to prove the other foliation perturbs to a tight contact structure, 
but this is much harder.
\ehw
It is very common for contact structures to be perturbations, even deformations, of foliations with Reeb components. In particular the following is known.
\begin{thm}[Etnyre 2006, \cite{EtnyreCFF}]
Every positive and negative contact structure on a closed oriented 3--manifold 
is a $C^\infty$-deformation of a $C^\infty$-foliation. 
Moreover, the foliation has Reeb components. 
\end{thm}
It would be very interesting to determine which contact structures were perturbations of Reebless foliations or taut (see below) foliations. 

A foliations $\xi$ on $M$ is called \dfn{taut} if each leaf of $\xi$ intersects a closed transversal curve. Equivalently, a foliation is taut if there is a vector field
$\xi$ that is transverse to $\xi$ and that preserves a volume form $\Omega$ on $M.$ From the first definition it is easy to see that a taut foliation is Reebless since
the torus leaf in a Reeb component separates $M$ and a transverse curve that passes through this leaf cannot close up (it's stuck on one side of the surface).
More generally, a taut foliation cannot contain a compact separating leaf. With this observation it is not hard to find Reebless foliations that are not taut. So tautness
is a strictly stronger notion that Reebless. 
\bhw
Try to prove the two definitions of taut are equivalent. \\ Hint: Poincar\'e recurrence is useful to prove the second definition implies the first. 
\ehw
\begin{thm}[Eliashberg-Thurston 1998, \cite{EliashbergThurston98}]\label{makefill}
If $\xi'$ is a positive contact structure that is $C^0$-close to a taut foliation $\xi$ then $\xi'$ is weakly symplectically semi-fillable. 
\end{thm}

\begin{proof}
Let $X=M\times [-1,1].$ Use $t$ as the coordinate on $[-1,1].$ Let $\alpha$ be a 1-form such that $\xi=\ker \alpha$ and set $\widetilde{\omega}=
\iota_v \Omega$ where $v$ and $\Omega$ are the vector field and volume form from the definition of taut. Note that $\widetilde{\omega}|_\xi>0$ and 
\[
d\widetilde{\omega}=d\iota_v\Omega= d\iota_v\Omega +\iota_v d\Omega=\mathcal{L}_v \Omega= 0\] 
(where $\mathcal{L}$ is the Lie derivative). Thus
$\omega=\widetilde{\omega}+\epsilon d(t\alpha)$ is a symplectic form on $X.$ 
\bhw
Show that $\omega\wedge \omega>0$ on $X.$
\ehw
Note also that $\omega|_{\xi\times\{\pm 1\}}>0.$ Thus if $\xi'$ is $C^0$-close to $\xi$ then $\omega|_{\xi'\times\{1\}}>0.$ Let $\xi''$ be a negative contact
structure on $M$ that is $C^0$ close to $\xi.$ Clearly $\omega|_{\xi''\times \{-1\}}>0$ so $(X,\omega)$ weakly symplectically fills  $(M, \xi')\coprod (-M,\xi'').$
\end{proof}

We can now construct lots of tight contact structures using a theorem of Gabai. 
\begin{thm}[Gabai 1983, \cite{Gabai83}] \label{GabaiThm}
Let $M$ be an irreducible 3-manifold and $\Sigma$ an oriented surface realizing a non-trivial homology class in $M$ and of minimal genus among 
representatives of its homology class. Then there is a taut foliation $\xi$ on $M$ with $\Sigma$ as a leaf.
\end{thm}
\begin{cor}
With $M$ and $\Sigma$ as above, there is a fillable contact structure $\xi'$ on $M$ such that $\langle e(\xi'),[\Sigma]\rangle=\pm \chi(\Sigma).$
\end{cor}
Actually, one needs to be a little careful proving this corollary. The foliation from Gabai's theorem is $C^2$ if the genus of $\Sigma$ is larger than 1, but when
the genus is 1 the foliation is not $C^2.$ However the only part of the foliation that is not $C^2$ is along the surface $\Sigma$ and there is holonomy along 
$\Sigma$. Thus looking back at the proof of Theorem~\ref{thm:f2c} we see we can still perturb this foliation into a contact structure. 

\section{Symplectic handle attachment and Legendrian surgery}\label{shandles}

We now want to discuss how to build symplectic manifolds. We will do this by starting with simple symplectic pieces and ``gluing'' them together. In order to 
accomplish this gluing we need a stronger notion of symplectic filling. To this end, let $(X,\omega)$ be a symplectic manifold, then we call a vector field 
$v$ \dfn{symplectically dilating} if 
\[
\mathcal{L}_v \omega=\omega.
\]
Suppose $v$ is transverse to $M=\partial X$ and is pointing out of $X$ along $M,$ then set $\alpha=(\iota_v \omega)|_M.$ We compute
\[
d\alpha= d\iota_v\omega=d\iota_v \omega+\iota_v d\omega=\mathcal{L}_v\omega=\omega,
\]
so
\[
\alpha\wedge d\alpha= (\iota_v \omega)\wedge \omega= \frac 12 \iota_v(\omega\wedge\omega).
\]
Since $\omega\wedge\omega$ is a volume form on $X$ it is clear that $\alpha\wedge d\alpha$ is a volume form on $M.$  In other words, $\alpha$ is a
contact form on $M.$ We say a contact manifold $(M,\xi)$ is \dfn{strongly filled} by a compact symplectic manifold $(X,\omega)$ if $\partial X=M$ and 
there is a dilating vector field $v$ for $\omega$ defined near the boundary of $X$ that is transversely pointing out of the boundary of $X$ and such that
$\iota_v\omega$ is a contact form for $\xi.$ We also say that $(X,\omega)$ is a \dfn{strong convex filling} of $(M,\xi).$ If the vector field $v$ points into
$X$ then we say $(X,\omega)$ is a \dfn{strong concave filling} of $(M,\xi).$
\bhw
Prove that a strong symplectic filling of a contact manifold $(M,\xi)$ is also a weak symplectic filling.  
\ehw
It is not true that a weak filling of a contact manifold is a strong filling, but we can sometimes create a strong filling from a weak filling.
\begin{thm}\label{weak2strong}
If $M$ is a homology sphere and $(X,\omega)$ is a weak symplectic filling of a contact structure $\xi$ on $M$ then $\omega$ can be altered to $\omega'$
so that $(X,\omega')$ is a strong symplectic filling of $(M,\xi).$
\end{thm}
\bhw
Try to prove this theorem. 
\ehw

\begin{figure}[ht]
  \relabelbox \small {\centerline{\epsfbox{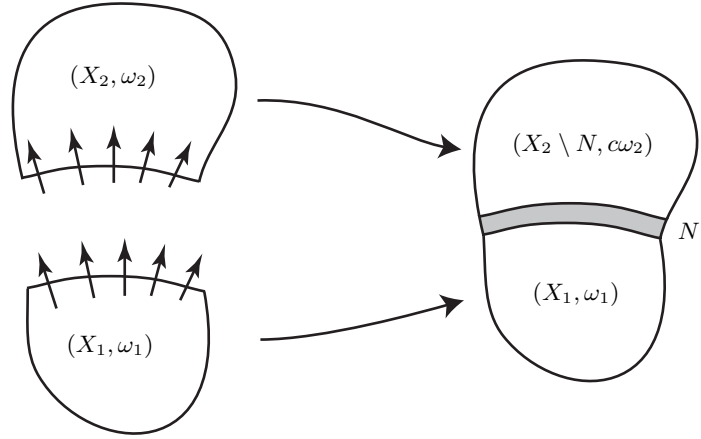}}} 
  \relabel{1}{$(X_1,\omega_1)$}
  \relabel{2}{$(X_2,\omega_2)$}
  \relabel{3}{$(X_2\setminus N,c\omega_2)$}
  \relabel{4}{$(X_1,\omega_1)$}
  \relabel{n}{$N$}  
  \endrelabelbox
        \caption{Gluing two symplectic manifolds.}
        \label{fig:glue}
\end{figure}

The key reason we have brought up the notion of strong fillings is to ``glue'' symplectic manifolds together. Specifically we have the following theorem. 
\begin{thm}\label{gluing}
If $(X_1,\omega_1)$ is a strong symplectic filling of $(M,\xi)$ and $(X_2,\omega_2)$ is a strong
concave filling of $(M,\xi)$ then $X=X_1\cup X_2$ has a symplectic structure $\omega$ such that
$\omega|_{X_1}=\omega_1$ and $\omega|_{X_2\setminus N}=c\omega_2$ where $N$ is a neighborhood of $\partial X_2$
in $X_2$ and $c>0$ is a
constant. See Figure~\ref{fig:glue}.
\end{thm}

Thus, after possibly rescaling the symplectic form on one of the pieces, a symplectic manifold with strongly convex boundary and one with strongly concave boundary can be glued together to get a symplectic manifold if the contact structures induced on their boundaries are the same. Such a gluing result is definitely not true for weak concavity/covexity. For more on this and all the various forms of convexity see \cite{Etnyre98}.

\begin{figure}[ht]
  \relabelbox \small {\centerline{\epsfbox{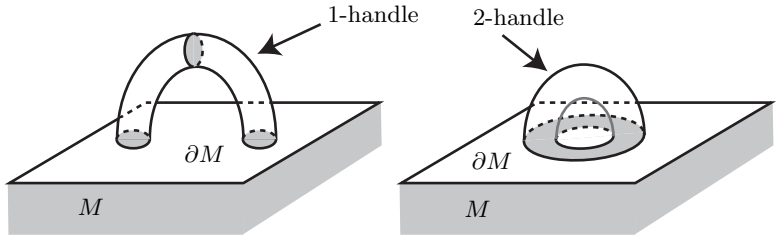}}} 
  \relabel{m}{$M$}
  \relabel{m1}{$M$}
  \relabel{b}{$\partial M$}
  \relabel{b1}{$\partial M$}
  \relabel{1}{$1$-handle}
  \relabel{2}{$2$-handle}  
  \endrelabelbox
        \caption{A 3-dimensional 1-handle, left, and 2-handle, right.}
        \label{fig:handle}
\end{figure}

We will use this gluing theorem repeatedly. Our first application of it will be to symplectic handle attachment. To this end we recall a little 4-dimensional topology. Suppose we are given a 4-manifold
$X$ with boundary. Then a \dfn{1-handle} is $h^1=D^1\times D^3$ (were $D^n$ is the unit disk in 
$\R^n$) and when we \dfn{attach a 1-handle} to $X$ we glue $h^1$ to $X$ along $A^1=(\partial D^1)\times D^3=S^0\times D^3=\{2 \text{ points}\}\times D^3.$ The set $A^1$ is called the attaching region 
of $h^1.$ If we identify two disjoint 3-balls in $\partial X$ there will be a unique way to 
glue $h^1$ to $\partial X$ so that $A$ goes to these two 3-balls. Since we cannot draw this picture we illustrate it one dimension lower in Figure~\ref{fig:handle}.

A \dfn{2-handle} is $h^2=D^2\times D^2$ and when we \dfn{attach a 2-handle} to $D$ we glue $h^2$ to $\partial X$ along $A^2=(\partial D^2)\times D^2=S^1\times D^2.$ Thus to attach a 2-handle we must
identify a knot $K$ in $\partial X,$ this is where $S^1\times \{pt\}$ will be glued, and we must also fix a framing on $K$ so that we will know how to glue $A^2$ to a neighborhood of $K.$ Recall, there is an integers worth of trivializations of $D^2\times S^1$ and any such trivialization is called a framing of the 
core of the solid torus.
\bhw
Show that an isotopy class of curves in $\partial (D^2\times S^1)$ that represents a generator in $\pi_1(D^2\times S^1)$ determines and is determined by a framing. 
\ehw
Again we illustrate, in Figure~\ref{fig:handle}, a 2-handle attachment in dimension three (here there is a 
unique way to attach a 2-handle once a curve in the boundary 3-manifold is specified).

Suppose $X'$ is $X$ with a 2-handle attached along a knot $K$ with framing $\mathcal{F}.$ Then
\[
\partial X'=[\partial X\setminus (S^1\times D^2)]\cup (D^2\times S^1)
\]
where the $S^1\times D^2$ is a neighborhood of $K.$ It is removed form $\partial X$ because when the
handle is attached it becomes part of the interior of $X'.$ Thus $\partial X'$ is obtained from $\partial X$
by removing a neighborhood of a knot $K$ and replacing it with another solid torus. This is clearly just 
a Dehn surgery on $\partial X.$ (If you have not studied Dehn surgery, then take this for the definition.)
\bhw
Check that this Dehn surgery is, in fact, a surgery with framing $\mathcal{F}.$ By this we mean that 
the new solid torus is glued in so that $\partial D^2\times\{pt\}$ is glued to a curve determined by $\mathcal{F}.$
\ehw

Now to relate this to contact geometry we say a knot $K$ in a contact manifold $(M,\xi)$ is \dfn{Legendrian} if $K$ is always tangent to $\xi$:
\[
T_xK\subset \xi_x,\quad x\in K.
\] 
Note that since $K$ is tangent to $\xi,$ the contact structure defines a trivialization of the normal bundle 
of $K$. That is, there is a natural \dfn{contact framing} on a Legendrian knot. 
\bhw
Convince yourself of this.\\
Hint: Take a vector field in $\xi$ along $K$ and use it to push off a copy of $K.$
\ehw
We can now discuss symplectic handle attachment.
\begin{thm}[Eliashberg, 1990 \cite{Eliashberg90a} and Weinstein, 1991 \cite{Weinstein91}]\label{thm:ew}
If $(X,\omega)$ is a symplectic manifold with strongly/weakly convex boundary and $X'$ is obtained
from $X$ by attaching a 1-handle to $X$ or attaching a 2-handle to $X$ along a Legendrian knot in $\partial X$ with framing one less than the contact framing, then $\omega$ extends to a symplectic
form $\omega'$ on $X'$ in such a way that $\partial X'$ has a strongly/weakly convex boundary.
\end{thm}
If $X'$ is obtained from $X$ by a symplectic 2-handle attachment along a Legendrian knot $K$ as in the theorem, 
then we say that the contact manifold $\partial X'$ is obtained from $\partial X$ by \dfn{Legendrian surgery} along $K.$
\begin{proof}[Sketchy proof of the theorem]
We consider 1-handle attachment.
The basic idea is to consider the model 1-handle in $\C^2$ illustrated in Figure~\ref{fig:shandle}. 
\begin{figure}[ht]
  \relabelbox \small {\centerline{\epsfbox{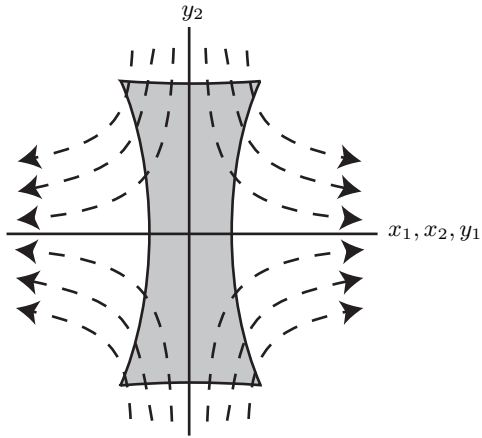}}} 
  \relabel{1}{$y_2$}
  \relabel{2}{$x_1,x_2,y_1$}
    \endrelabelbox
        \caption{A model symplectic 1-handle.}
        \label{fig:shandle}
\end{figure}
Here we give
$\C^2$ coordinates $(z_1=x_1+iy_1, z_2=x_2+iy_2).$ In this model we can construct an expanding vector field
$v'$ that is transversely pointing into the handle along the attaching region $A^1$ and transversely 
pointing out of the handle along the other boundary component. Now if $v$ is the vector field implicated
in the definition of the convexity of $\partial X$ then we can use $v'$ and $v$ to glue the standard symplectic structure on $\C^2$ to $\omega$ when we attach $h^1$ to $X.$ You can do this by 
carrying out the following exercises.
\bhw
Show that the contact structures induced on the attaching region $A^1\subset h^1$ and on the two 3-balls
to which $A^1$ is glued are the same. Moreover, show that we can choose our gluing map so that the contact forms are the same.
\ehw
\bhw
Use the previous exercise and the flow of $v$ and $v'$ to show that a neighborhood of $A^1$ in 
$\C^2\setminus h^1$ is symplectomorphic to a neighborhood of the two 3-balls in $X.$ 
Use this symplectomorphism to glue $h^1$ to $X$ and extend $\omega$ over the handle.
\ehw

Gluing a 2-handle is a little more tricky, but the main idea is the same.
\bhw
Try to extend $\omega$ over a 2-handle attached as in the statement of the theorem. Or at least try to 
figure out why you have the framing condition in the theorem. 
\ehw

The observant reader will have noticed that we seem to be assuming that $X$ has strongly convex
boundary, but the theorem works for weakly convex boundary as well. This follows from (a relative version of) Theorem~\ref{weak2strong}.
\end{proof}

\section{Open book decompositions and symplectic caps}\label{sec:ob}
Our main goal in this section is to sketch the proof of the following theorem. 
\begin{thm}[Eliashberg and Etnyre, 2004 \cite{Eliashberg04, Etnyre04a}]\label{eet}
If $(X,\omega)$ is a compact symplectic manifold with weakly convex boundary then there is a closed symplectic manifold $(X',\omega')$ into
which $(X,\omega)$ symplectically embeds. 
\end{thm}
Recall from the Introduction that this completes the  contact geometric input into the various advance in low-dimensional topology discussed there.
We will only sketch the ideas in the proof of this theorem here, for a more complete discussion see the original papers \cite{Eliashberg04, Etnyre04a} 
or the survey article \cite{EtnyreOBN}.

The last ingredient we need to prove this theorem is open book decompositions. Let $\Sigma$ be a compact oriented surface with boundary and
$\phi:\Sigma\to \Sigma$ an orientation preserving diffeomorphism of $\Sigma$ that is equal to the identity near $\partial \Sigma.$ The 
\dfn{mapping torus} of $\phi$ is 
\[T_\phi=\Sigma\times[0,1]/\sim,\]
where $(x,1)\sim (\phi(x),0).$ For each boundary component of $\Sigma$ we have a boundary component of $T_\phi.$ This boundary component
is a torus, with two canonical circles $\mu=\{pt\}\times [0,1]/\sim$ and $\lambda=\partial \Sigma\times \{pt\}.$ Let 
\[
M_{(\Sigma,\phi)}=\left(T_\phi\cup \coprod_{|\partial \Sigma|} S^1\times D^2\right)/\sim,
\]
where $\sim$ glues a solid torus $S^1\times D^2$ to $T_\phi$ so that $\{pt\}\times \partial D^2$ is glued to $\mu$ in $\partial T_\phi$ and $S^1\times \{pt\}$
is glued to $\lambda.$
\bhw
If $L$ is the union of the cores of all the $S^1\times D^2$'s in $M_{(\Sigma,\phi)}$ then show that $M\setminus L$ fibers over the circle with fiber diffeomorphic
to $\Sigma.$ We call $L$ the \dfn{binding} of the open book and $\Sigma$ the \dfn{page}.
\ehw

An \dfn{open book decomposition}, or just open book, of a closed oriented 3-manifold $M$ is an identification of $M$ with $M_{(\Sigma,\phi)}$ for some $(\Sigma,\phi)$
as above. This is not the best definition of open book as it only defines an open book up to diffeomorphism, but it will suffice for our purposes. For
a better definition see \cite{EtnyreOBN}.

\begin{fct} All closed oriented 3-manifolds have open book decompositions 
\end{fct}
\bhw
Prove this fact!\\
Hint: You might find it useful to recall that all 3-manifold are branched covers of $S^3$ branched over some link and that links can be braided about
the unknot.
\ehw

An open book $(\Sigma, \phi)$ for $M$ is said to \dfn{support} or be \dfn{compatible with} a contact structure $\xi$ on $M$ if there is a 1-form $\alpha$ such that
$\xi=\ker\alpha,$ $\alpha (v)>0$ for any $v\in TL$ that agrees with the orientation on the binding $L,$ and $d\alpha|_{page}\not=0$ and induces the correct
orientation on the page.  Thurston and Winkelnkemper \cite{ThurstonWinkelnkemper75} have shown that every open book supports a contact structure and Giroux \cite{Giroux02} observed this
contact structure is unique.

Given an open book $(\Sigma,\phi)$ supporting a contact structure $\xi$ on $M$ then the \dfn{positive stabilization of $(\Sigma,\phi)$} is the open book $(\Sigma',
\phi')$ with
\[
\Sigma'=\Sigma \cup (\text{1-handle})
\]
and $\phi'=\phi\circ D_\gamma$ where $D_\gamma$ is a right handed Dehn twist along a curve $\gamma$ that runs over the new 1-handle in $\Sigma'$ exactly once. 
\bhw
Show that $(\Sigma',\phi')$ is still an open book for $M$ and still supports $\xi.$
\ehw
\begin{thm}[Giroux 2002, \cite{Giroux02}]\label{GirouxCor}
There is a one-to-one correspondence between 
\[
\{\text{oriented contact structures on $M$ up to isotopy}  \}
\]
{\centerline {and}}
\[
\{\text{open book decompositions of $M$ up to positive stabilization} \}.
\]
\end{thm}

To use this theorem to prove Theorem~\ref{eet} we need to see how Legendrian surgery interacts with open books. To this end suppose $(\Sigma,\phi)$ supports 
$\xi$ on $M.$ We begin by forgetting  about the contact structure $\xi$ and concentrating on $M.$ Let $\gamma$ be a simple closed curve contained
in a page of the open book. Note that $\gamma$ gets a framing $\mathcal{F}$ from the page. Let $M'$ be the manifold obtained from $M$ by $\mathcal{F}\pm 1$
surgery on $\gamma.$ 
\bhw
Show that an open book for $M'$ is $(\Sigma, \phi\circ D^\mp_\gamma).$\\
Hint: If you cut $M\setminus L$ open along the page containing $\gamma$ and reglue by $D^\pm_\gamma$ then the resulting manifold will differ from $M$
in a neighborhood of $\gamma.$ The key now is to see that this difference is a Dehn surgery with the appropriate framing. 
\ehw
\begin{fct}[Legendrian realization principle \cite{Honda00a}]
If $\gamma$ is a non-separating curve on a page of the open book then we can isotop the open book slightly so that $\gamma$ is Legendrian and the 
contact framing agrees with the page framing.
\end{fct}
Actually this is not exactly the Legendrian realization principle, but this fact easily follows from it.

\begin{fct}\label{changebylsurg}
Let $\gamma$ be a Legendrian knot in a page of the open book $(\Sigma,\phi)$ supporting the contact manifold $(M,\xi).$ If $(M',\xi')$ is obtained from $(M,\xi)$ by Legendrian surgery on $\gamma,$ then  
$(M',\xi')$ is supported by $(\Sigma,\phi\circ D_\gamma).$
\end{fct}

So, summarizing,  given a symplectic filling $(X,\omega)$ of $(M,\xi)$ and a non-separating curve $\gamma$ on a page of an open book $(\Sigma,\phi)$ supporting $\xi$ then we
can attach a symplectic 2-handle to $(X,\omega),$ as in Theorem~\ref{thm:ew}, to get $(X',\omega')$ and $\partial X'=M'$ and $\xi'$ is filled by $\omega'$ where a
supporting open book for $\xi'$ is $(\Sigma,\phi\circ D_\gamma).$ 

Now we need some facts about the mapping class group of a surface. 
\begin{fct}\label{fct:mapping}
If $\Sigma$ is a surface with one boundary component then any diffeomorphism of $\Sigma$ that is fixed on the boundary can be written as
\[\phi=D_c^m\circ D_{\gamma_1}^{-1}\circ\cdots\circ D_{\gamma_n}^{-1},\]
where the $\gamma_i$ are separating curves on $\Sigma$ and $c$ is a curve parallel to $\partial \Sigma.$
\end{fct}

We are now ready to begin our sketch of the proof of Theorem~\ref{eet}. We start with a weak symplectic filling $(X,\omega)$ of the contact manifold
$(M,\xi).$ Let $(\Sigma,\phi)$ be an open book for $M$ supporting $\xi.$ By positively stabilizing if necessary we can assume that $\Sigma$ has only
one boundary component. We can further assume
\[
\phi=D_c^m\circ D_{\gamma_1}^{-1}\circ\cdots\circ D_{\gamma_n}^{-1},
\]
as in Fact~\ref{fct:mapping} above.

\begin{figure}[ht]
  \relabelbox \small {\centerline{\epsfbox{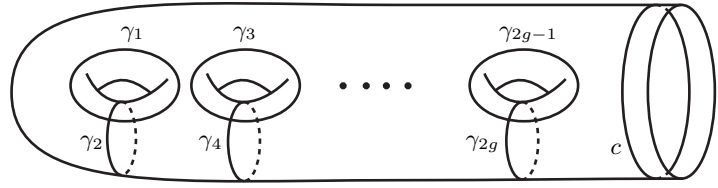}}} 
  \relabel{1}{$\gamma_1$}
  \relabel{2}{$\gamma_2$}
   \relabel{3}{$\gamma_3$}
  \relabel{4}{$\gamma_4$}
   \relabel{5}{$\gamma_{2g-1}$}
  \relabel{6}{$\gamma_{2g}$}
  \relabel{7}{$c$}
    \endrelabelbox
        \caption{The curves $\gamma_i$ and $c$ on $\Sigma.$}
        \label{fig:sfc1}
\end{figure}

\begin{figure}[ht]
  \relabelbox \small {
  \centerline{\epsfbox{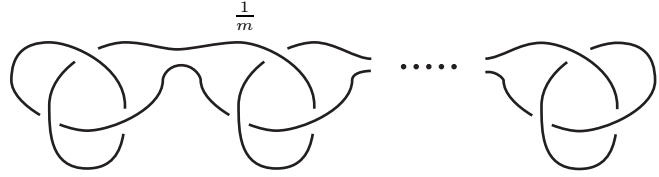}}}
  \relabel{m}{$\frac{1}{m}$} \endrelabelbox
        \caption{Topological description of $M''.$}
        \label{fig:spic}
\end{figure}

We would now like to simplify the monodromy of this open book. Using Theorem~\ref{thm:ew} we can attach 2-handles to $X$ along Legendrian knots in the pages of the open book for $(M,\xi)$ 
and extend the symplectic structure of $X$ over the 2-handles. Then since the upper boundary of the new 
4-manifold is obtained from $(M,\xi)$ by Legendrian surgery Fact~\ref{changebylsurg}  tells us that the monodromy of this new boundary will have an extra right handed Dehn twists along each of the attaching curves of the 2-handles. Thus if we can attach $n$ symplectic 2-handles to $X$ along the $\gamma_i$'s to get a symplectic manifold $(X',\omega')$ with $\partial (X',\omega')=
(M',\xi')$  and $\xi'$ supported by the open book $(\Sigma,\phi')$ where 
\[\phi'=D^m_c.\]
Now attach $2g$ more symplectic 2-handles to get $(X'',\omega'')$ having contact boundary $(M'',\xi'')$ supported by the open book $(\Sigma,\phi'')$ where
\[\phi''=D^m_c\circ D_{\gamma_1}\circ\cdots\circ D_{\gamma_{2g}},\]
where the curves $\gamma_i$ are shown in Figure~\ref{fig:sfc1}.
\bhw
Show that if $\Sigma$ has genus $g$ then $M''$ is the manifold depicted in Figure~\ref{fig:spic}. In particular, $M''$ is a homology 3-sphere. 
\ehw
Thus according to Theorem~\ref{weak2strong} we can slightly perturb $\omega''$ so that  $(X'',\omega'')$ is a strong symplectic filling of $(M'',\xi'').$ Also note
that since $(X'',\omega'')$ was built by attaching symplectic 2-handles to $(X,\omega),$ we can symplectically embed $(X,\omega)$ into $(X'',\omega'').$
We can now hope to glue a concave filling of $(M'',\xi'')$ to $(X'',\omega'')$ to get a closed symplectic manifold. Before we do this, we will attach a few more symplectic
2-handles to $(X'',\omega'')$ to further normalize its boundary. 

Recall if  $\gamma_1,\ldots, \gamma_{2g}$ is a collection of  simple closed curves in $\Sigma$ that satisfy $\gamma_i\cdot \gamma_j$ is $1$ 
if $|i-j|=1$ and is 0 otherwise, where $\cdot$ means
geometric intersection, then $\Sigma$ is a neighborhood of the union of the $\gamma_i$'s. Moreover we have the so called chain relation in the mapping class
group of $\Sigma$:
\[
(D_{\gamma_1}\circ \ldots \circ D_{\gamma_{2g}})^{4g+2} = D_c.
\]
Using the chain relation it is easy to see that we can attach $8g^2+3g$ more symplectic 2-handles to $(X'',\omega'')$ to obtain $(X^{(3)},\omega^{(3)})$ which strongly
fills $(M^{(3)},\xi^{(3)})$ where $\xi^{(3)}$ is supported by $(\Sigma,\phi^{(3)})$ where 
\[
\phi^{(3)}=D_c^{k},
\]
for some $k.$

Finally, we can stabilize the open book $(\Sigma,\phi^{(3)})$ for $(M^{(3)},\xi^{(3)})$ so that we can attach more symplectic 2-hanldes to $(X^{(3)}\omega^{(3)})$ to get a symplectic
manifold $(X^{(4)},\omega^{(4)})$ strongly symplectically filling $(M^{(4)},\xi^{(4)})$ where $\xi^{(4)}$ is supported by $(\Sigma',\phi^{(4)})$ where 
\[
\phi^{(4)}=D_{c'},
\]
and $c'$ is again a curve parallel to the boundary of $\Sigma'.$
\bhw
Prove this last assertion.
\ehw
\bhw
Show that $M^{(4)}$ is an $S^1$-bundle over $\Sigma''$ with Euler number $-1$ where $\Sigma''$ is the closed surface obtained from  $\Sigma'$ by attacing 
a disk to $\partial \Sigma'.$
\ehw
\bhw
If $Y$ is the $D^2$ bundle over $\Sigma''$ with Euler number 1, then $Y$ admits a symplectic structure with concave boundary and $\partial Y=-M^{(4)}.$\\
Hint: Think about a connection on the circle bundle $\partial Y.$
\ehw
\bhw
Show the contact structure induced on $\partial Y$ from the strong concave filling is contactomorphic to $\xi^{(4)}.$\\
Hint: Both contact structures are transverse to the fibers of the $S^1$-bundle.
\ehw

Now we can use Theorem~\ref{thm:ew} to 
glue $Y$ with the above constructed symplectic structure to $(X^{(4)},\omega^{(4)})$ to get a closed symplectic manifold into which $(X,\omega)$ embeds, thus
completing the proof of Theorem~\ref{eet}.

\def\cprime{$'$} \def\cprime{$'$}

\end{document}